\documentclass[12pt]{article}

\usepackage[utf8]{inputenc}
\usepackage[T1]{fontenc}
\usepackage{palatino}
\usepackage{amsmath}
\usepackage{amssymb}
\usepackage{embedfile}
\usepackage{booktabs}
\usepackage{amsthm}
\usepackage{thmtools}
\usepackage{abstract}
\usepackage{tikz-cd}
\usepackage{enumerate}
\usepackage{hyperref}
\usepackage[capitalize]{cleveref}
\usepackage[top=1.2in]{geometry}
\usepackage{mathtools}

\crefformat{equation}{Equation~\textup{#2#1#3}}
\Crefformat{equation}{Equation~\textup{#2#1#3}}
\crefrangeformat{equation}{Equations~\textup{#3#1#4} to~\textup{#5#2#6}}
\Crefrangeformat{equation}{Equations~\textup{#3#1#4} to~\textup{#5#2#6}}

\usepackage[font=footnotesize]{caption}

\numberwithin{equation}{subsection}

\declaretheoremstyle[
    spacebelow=2\topsep,
    spaceabove=2\topsep,
    headfont=\normalfont\bfseries,
    bodyfont=\itshape,
    postheadspace=\newline,
    qed=${\lrcorner}$,
    headpunct={},
    notebraces={[}{]}
]{breakit}

\declaretheoremstyle[
    spacebelow=2\topsep,
    spaceabove=2\topsep,
    headfont=\normalfont\bfseries,
    bodyfont=\normalfont,
    postheadspace=\newline,
    qed=${\lrcorner}$,
    headpunct={},
    notebraces={[}{]}
]{breakup}

\declaretheorem[numberlike=equation,style=breakit]{theorem}
\declaretheorem[numberlike=theorem,style=breakit]{lemma}

\declaretheorem[numberlike=theorem,style=breakup]{definition}
\declaretheorem[numberlike=theorem,style=breakup]{example}
\declaretheorem[numberlike=theorem,style=breakup]{note}

\numberwithin{theorem}{subsection}
\numberwithin{lemma}{subsection}
\numberwithin{corollary}{subsection}

\numberwithin{definition}{subsection}
\numberwithin{example}{subsection}
\numberwithin{note}{subsection}

\newcommand{\zz}{\mathbb{Z}}
\newcommand{\rr}{\mathbb{R}}

\newcommand{\nn}{\mathbb{N}}
\newcommand{\rp}{\mathbb{R}\mathbb{P}}

\DeclareMathOperator{\im}{im}
\newcommand{\rrleq}{\rr_\leqslant}
\newcommand{\Vecrrdiag}{[\rr_\leqslant,\mathsf{Vec}]}
\newcommand{\ccat}{\mathcal{C}}
\newcommand{\dcat}{\mathcal{D}}
\newcommand{\Hom}{\mathrm{Hom}}
\newcommand{\Top}{\mathsf{Top}}
\DeclareMathOperator{\id}{id}
\DeclareMathOperator{\dist}{d}
\newcommand{\nt}{\Rightarrow}
\newcommand{\vare}{\varepsilon}
\newcommand{\Fup}{F_\uparrow}
\newcommand{\Fdown}{F^\downarrow_s}
\newcommand{\Fpaira}{(\Fup(a),\Fdown(a))}

\usepackage{needspace}
\usepackage{titlesec}

\titleformat*{\section}{\needspace{2in}\centering\Large\bfseries}
\titleformat*{\subsection}{\needspace{1in}\large\bfseries}
\titleformat*{\subsubsection}{\needspace{1in}\large\bfseries}

\usepackage{appendix}
\crefname{appsec}{Appendix}{Appendices}

\title{\vspace{-1em}Death and extended persistence in computational algebraic topology}
\author{Timothy Hosgood}

\begin{document}
    \embedfile{\jobname.tex}

    \maketitle

    \begin{abstract}
        \small
        The main aim of this paper is to explore the ideas of \emph{persistent homology} and \emph{extended persistent homology}, and their \emph{stability theorems}, using ideas from \cite{CohenSteiner:2009ho,CohenSteiner:2007is,Bubenik:dn}, as well as other sources.
        The secondary aim is to explore the homology (and cohomology) of non-orientable surfaces, using the Klein bottle as an example.
        We also use the Klein bottle as an example for the computation of (extended) persistent homology, referring to it throughout the paper.
    \end{abstract}
    
    \renewcommand{\abstractname}{Contents}
    \begin{abstract}
        \noindent\vspace{-2em}\tableofcontents
    \end{abstract}

    \section{Introduction} 
    \label{sec:introduction}

        \begin{quotation}
            \raggedleft
            \emph{``Our birth is nothing but our death begun.''}\\
            -- Edward Young, \emph{Night Thoughts}
        \end{quotation}

        The main aim of this paper is to explore the ideas of \emph{persistent homology} and \emph{extended persistent homology}, and their \emph{stability theorems}, using ideas from \cite{CohenSteiner:2009ho,CohenSteiner:2007is,Bubenik:dn}, as well as other sources.
        The secondary aim is to explore the homology (and cohomology) of non-orientable surfaces, using the Klein bottle as an example.
        We also use the Klein bottle as an example for the computation of (extended) persistent homology, referring to it throughout the paper.
        There are numerous diagrams and sketches, as well as small computational examples, in the hope that the topological nature of this subject doesn't get lost amidst the algebra.

        A lot of consideration has been given to ensuring that this paper is as self-contained as possible (without being overly long) but whilst mentioning other (recent) papers, since many of the ideas found in this paper are relatively modern.
        In particular, and as has always been the case with algebraic topology, the subject is leaning more and more towards category-theoretic language -- many ideas that haven't been around for very long are already being rephrased in new ways.
        We try to place equal emphasis on both approaches, drawing inspiration from \cite{CohenSteiner:2009ho} for the topological view, and \cite{Bubenik:dn} for the category-theoretic view.
        In a sense, this paper aims to be an addendum to \cite{Edelsbrunner:2008gf}, which is a brilliant survey of persistent homology, in light of some of the results from \cite{Bubenik:dn}.

        \subsection{Conventions and notation} 
        \label{sub:conventions}

            Unless otherwise stated, we adopt the following conventions and notation:
            \begin{itemize}
                \item all (co)homology\footnote{
                    We write `(co)homology' to mean `homology and cohomology'.
                } groups have coefficients in $\zz/2\zz$;
                \item for a group $G$ we write $G^n$ to mean $\bigoplus_{i=1}^n G$;
                \item $\overline{D^n}=\{x\in\rr^n : \|x\|\leqslant1\}$ is the \emph{closed} $n$-ball or $n$-disc;
                \item $D^n=\{x\in\rr^n : \|x\| < 1\}=(\overline{D^n})^\circ$ is the \emph{open} $n$-ball or $n$-disc;
                \item $S^{n-1}=\{x\in\rr^n : \|x\|=1\}=\partial(\overline{D^n})$ is the $n$-sphere;
                \item we (\emph{sometimes}\footnote{
                    This is not a strict convention, but we often use this shorthand to save space.
                }) write $\zz_n$ to mean $\zz/n\zz$;
                \item we write $G\langle x_1,x_2,\ldots,x_n\rangle$ to mean the group $G^n$ with basis $x_1,x_2,\ldots,x_n$;
                \item for a category $\ccat$ we write $x\in\ccat$ to mean that $x\in\mathrm{ob}(\ccat)$ is an object of $\ccat$;
                \item when we say `an interval $I\subseteq\rr$' we mean any interval, i.e. open, closed, half-open half-closed, or even infinite;
                \item if $\eta\colon F\nt G$ is a natural transformation between functors then we write $\eta(x)$ to mean the constituent morphism $\eta(x)\colon F(x)\to G(x)$;
                \item if we write $a\geqslant0$ then we mean, in particular, that $a\in\rr$;
                \item we write $\{*\}$ to mean a singleton (a set with one element);
                \item $0\in\nn$.
            \end{itemize}


        \subsection{Background knowledge} 
        \label{sub:assumptions}

            We assume that the reader has a knowledge of some of the fundamental notions in algebraic topology, namely: simplicial- and $\Delta$-complexes, and simplicial and singular homology (both relative and absolute).
            All of these topics are covered in \cite[\S2.1]{hatcher2002algebraic}.
            In particular, we use the following theorem.

            \begin{theorem}\label{th:simplicial-iso-singular}
                Let $X$ be a $\Delta$-complex.
                Then the $k$-th simplicial homology group is isomorphic to the $k$-th singular homology group, i.e.
                \begin{equation*}
                    H_k(X)\cong H_k^\Delta(X).\qedhere
                \end{equation*}
            \end{theorem}

            \begin{proof}
                This is a specific case (where $A=\varnothing$) of \cite[Theorem~2.27,~\S2.1]{hatcher2002algebraic}.
            \end{proof}

            Because of this theorem, for any $\Delta$-complex $X$ we can write $H_k(X)$ to mean \emph{the $k$-th homology group of $X$} without specifying whether it is calculated using simplicial or singular homology -- up to isomorphism, the two are the same.



    \section{The (co)homology of non-orientable surfaces} 
    \label{sec:the_}

        \begin{quotation}
            \raggedleft
            \emph{``Which way is up, what's goin' down? I just don't know, no''}\\
            -- Barry White, \emph{Which Way Is Up}
        \end{quotation}

        Our first aim is to explore the (co)homology of non-orientable surfaces, using the Klein bottle as an explicit example.
        Two of the mains tools that we use are \emph{cellular homology} and \emph{Poincaré duality}; we summarise most of the necessary definitions and results, as well as defining some non-standard notation, in \cref{sec:appendix_cells}.

        \subsection{Non-orientable surfaces} 
        \label{sub:non_orientable_surfaces}

            \begin{definition}[Non-orientable surface of genus $g$]
                For $g>0$ let $N_g$ be the surface obtained from a regular $2g$-gon by identifying its edges according to a cyclic labelling of its edges $a_1a_1a_2a_2\ldots a_ga_g$.
            \end{definition}

            \begin{figure}[ht]
                \centering
                \includegraphics[width=.4\textwidth]{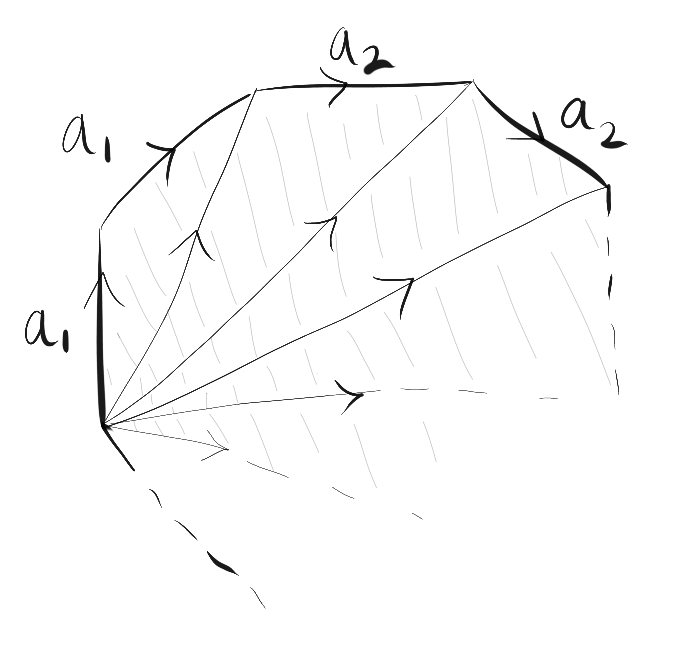}
                \caption{An sketch of a $\Delta$-complex structure on $N_g$, where the arrows indicate the ordering of the simplices. If we apply barycentric subdivision twice to this structure (or indeed, to any $\Delta$-complex) then we obtain a simplicial-complex structure \cite[Exercise~23,~\S2.1]{hatcher2002algebraic}.}\label{fg:ng-and-nd-delta}
            \end{figure}

            We know that we can also consider $N_g$ as a CW complex by using its construction as a polygon with pairwise side identification (see \cref{sub:cw_complexes_and_cellular_homology}).
            Explicitly, the CW-complex structure on $N_g$ has the following properties:
            \begin{enumerate}[(i)]
                \item $\dim_\text{CW}N_g=2$;
                \item $\sigma(N_g) = \{e_0^0\} \cup \{e_0^1,e_1^1,\ldots,e_{g-1}^1\} \cup \{e_0^2\}$;
                \item $\varphi_\alpha^1\colon S^0=\{-1,1\}\to N_g^0=\{e_0^0\}$ is the constant map;
                \item $\varphi_0^2\colon S^1\to N_g^1=\vee_{i=1}^g S^1$ is the map $a_1a_1a_2a_2\ldots a_ga_g$.
            \end{enumerate}
            We can use this to calculate the (co)homology of $N_g$ (recalling \cref{th:cellular-iso-sing-or-simp}).

            \bigskip

            By property (ii)\footnote{
                Recalling \cite[Lemma~2.34,~\S2.2]{hatcher2002algebraic}: $C_n^\text{CW}(X)$ is free abelian with basis in bijective correspondence to $n$-cells of $X$.
            }, the associated cellular chain complex of $N_g$ is
            \begin{equation}\label{eq:cell-chain-complex-ng}
                \begin{tikzcd}
                    0 \arrow{r} & \zz/2\zz \arrow{r}{d_2} & (\zz/2\zz)^g \arrow{r}{d_1} & \zz/2\zz \arrow{r} & 0.
                \end{tikzcd}
            \end{equation}
            We need to calculate the maps $d_2$ and $d_1$ to compute the homology groups of this chain complex.
            Since $N_g^0$ is a singleton set and $N_g$ is connected, we know\footnote{
                Because otherwise $H_0(X)$ would not be $\zz$.
            } that $d_1^\Delta=0$, and so \cref{th:cellular-boundary-formula} tells us that $d_1=0$.
            To calculate $d_2$ we just need to know $\deg(\chi^2_{0\beta})$ for all $\beta\in\{0,1,\ldots,g-1\}$, also by \cref{th:cellular-boundary-formula}.
            Here it is easiest to calculate the degree using the \emph{local degree}\footnote{
                See \cite[Proposition~2.30,~\S2.2]{hatcher2002algebraic} and the preceding paragraphs.
            }.
            For any point $p\in S^1$ we have two points in the preimage under $\chi_{0\beta}^2$.
            Since the attaching map just `wraps around' circles, we see that the attach-and-collapse map is of local degree 1.
            Putting these two facts together we see that $\deg(\chi_{0\beta}^2)=2=0$ in $\zz_2$, and so $d_2=0.$

            Because both $d_1$ and $d_2$ are zero, we can read the homology groups straight off from \cref{eq:cell-chain-complex-ng}; using \cref{th:poincare-duality} we can calculate the cohomology groups from the homology groups.
            This gives the following.
            \begin{equation}\label{eq:co-homology-of-ng}
                H_k(N_g) = H^k(N_g) =
                \begin{cases}
                    \,\,\zz/2\zz & k=0,2\\
                    (\zz/2\zz)^g & k=1\\
                    \,\,0 & k\geqslant3.
                \end{cases}
            \end{equation}


        \subsection{The (co)homology of the Klein bottle} 
        \label{sub:the_klein_bottle_co_homology}

\begin{figure}[ht]
                \centering
                \includegraphics[width=.5\textwidth]{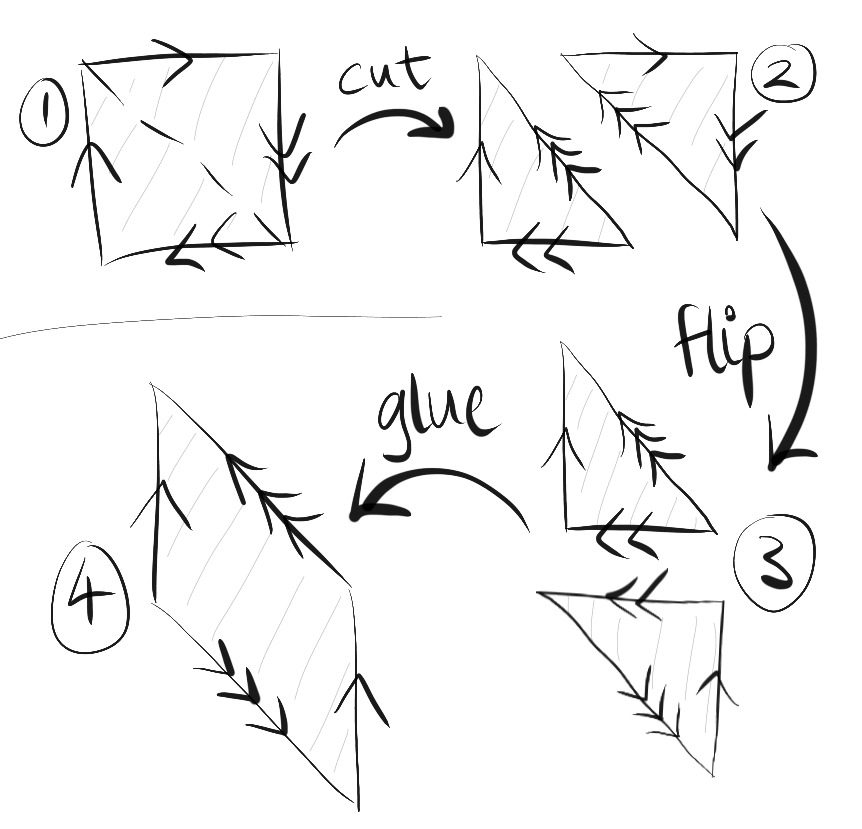}
                \caption{Showing that $N_2$ is homeomorphic to the Klein bottle $K$ by `cutting and gluing'. Here the arrows represent side identification.}\label{fg:n2-cong-k}
            \end{figure}

            We claim\footnote{
                More generally, $N_g$ is homeomorphic to the connected sum of $g$ copies of $\rp^2$.
                This can be seen from the fact that $\rp^2$ can be described as the 2-gon with boundary word $aa$, and that the connected sum of two polygons with boundary words $x_1\ldots x_m$ and $y_1\ldots y_n$ (respectively) is $x_1\ldots x_my_1\ldots y_n$.
                See \cite[\S1.5]{Massey:1967we}.
            } that $N_2$ is homeomorphic to $K$, where $K$ is the Klein bottle: the non-orientable surface constructed from a square with side identifications $abab^{-1}$ (see \cref{fg:n2-cong-k}).

\begin{figure}[ht]
                \centering
                \includegraphics[width=.3\textwidth]{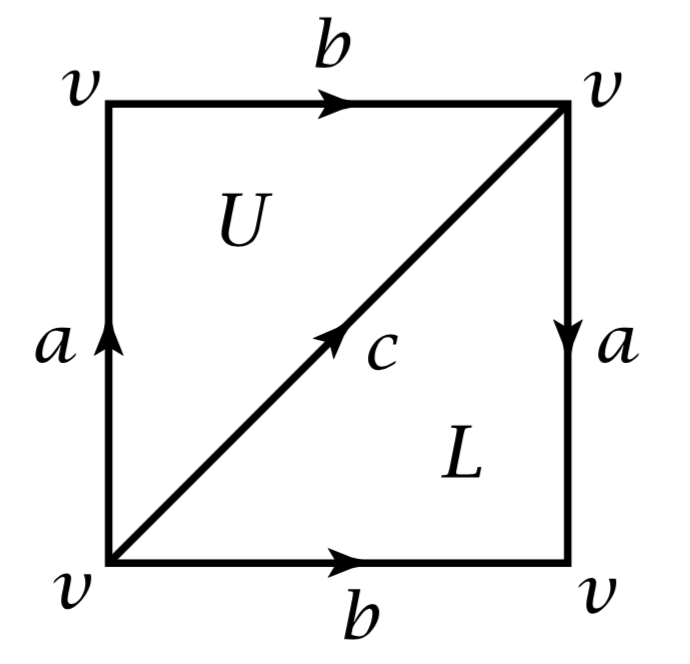}
                \caption{The Klein bottle $N_2$ with a $\Delta$-complex structure \cite[p.~102,~\S2.1]{hatcher2002algebraic}.}\label{fg:hatcher-klein-bottle}
            \end{figure}

            Using \cref{eq:co-homology-of-ng} with $g=2$, we see that the Klein bottle $N_2$ has (co)homology
            \begin{equation}
                H_k(N_2) = H^k(N_2) =
                \begin{cases}
                    \zz_2 & k=0,2\\
                    \zz_2\oplus\zz_2 & k=1\\
                    0 & k\geqslant3.
                \end{cases}
            \end{equation}
            However, if we wish to find explicit cycle representatives for the generators of (the non-trivial) $H_k(N_2)$ it is easier to use the $\Delta$-complex structure of $N_2$ (see \cref{fg:hatcher-klein-bottle}) and simplicial homology.

            We have the simplicial chain complex
            \begin{equation*}
                0\to\Delta_2(N_2)\xrightarrow{\partial_2}\Delta_1(N_2)\xrightarrow{\partial_1}\Delta_0(N_2)\to0
            \end{equation*}
            which (using the labelling from \cref{fg:hatcher-klein-bottle}) takes the explicit form
            \begin{equation}\label{eq:simplicial-chain-complex-klein-bottle}
                \begin{array}{cccccccccc}
                    0 &\to &\zz_2\langle U,L\rangle &\xrightarrow{\partial_2} &\zz_2\langle a,b,c\rangle &\xrightarrow{\partial_1} &\zz_2\langle v\rangle &\to &0.\\[0.3em]
                    & &U &\mapsto &a+b+c & & & &\\
                    & &L &\mapsto &a+b+c & & & &\\
                    & & & &a &\mapsto &2v=0 & &\\
                    & & & &b &\mapsto &0 & &\\
                    & & & &c &\mapsto &0 & &\\
                \end{array}
            \end{equation}
            Using \cref{eq:simplicial-chain-complex-klein-bottle} we can read off explicit generators for the homology groups:
            \begin{align}
                H_0(N_2) &= \zz_2\langle v\rangle;\\
                H_1(N_2) &= \zz_2\langle a,b\rangle;\label{eq:explicit-homology-1-klein-bottle}\\
                H_2(N_2) &= \zz_2\langle U+L\rangle,
            \end{align}
            where \cref{eq:explicit-homology-1-klein-bottle} comes from the fact that
            \begin{equation*}
                H_1(N_2) = \frac{\ker\partial_1}{\im\partial_2} = \frac{\zz_2\langle a,b,c\rangle}{\zz_2\langle a+b+c\rangle} \cong \frac{\zz_2\langle a,b,a+b+c\rangle}{\zz_2\langle a+b+c\rangle} \cong \zz_2\langle a,b\rangle,
            \end{equation*}
            since $(a+b+c)+a+b=c$.



    \section{Persistent homology} 
    \label{sec:persistent_homology}

        \begin{quotation}
            \raggedleft
            \emph{``The only thing wrong with immortality is that it tends to go on forever.''}\\
            --Herb Caen, \emph{Herb Caen's San Francisco}
        \end{quotation}

        \emph{Many of the definitions in this section come from \cite{Bubenik:dn}.}

        \bigskip

        We now consider the \emph{persistent homology} of topological spaces under certain (reasonably weak) hypotheses.
        We assume that the reader is familiar with some basic concepts of category theory, such as functors and natural transformations.

        \begin{definition}[Functor categories and diagrams]
            Let $\ccat,\dcat$ be categories where $\ccat$ is small\footnote{
                i.e. the objects of $\ccat$ form a set, not a proper class.
            }.
            Define $[\ccat,\dcat]$ to be the \emph{functor category}: its objects are functors $F\colon\ccat\to\dcat$, called \emph{$\ccat$-indexed diagrams in $\dcat$}, or \emph{\mbox{$[\ccat,\dcat]$ diagrams}}; its morphisms are natural transformations between such functors.
        \end{definition}

        \begin{note}
            Many sources (including \cite{Bubenik:dn}) use the notation $\dcat^\ccat$ instead of $[\ccat,\dcat]$.
            We use the latter simply as a matter of upbringing.
        \end{note}

        \begin{definition}[Poset categories]
            Let $(P,\leqslant)$ be a poset\footnote{
                i.e. a set $P$ equipped with a partial order $\leqslant$.
            }.
            Define the \emph{poset category $P_\leqslant$} as follows: its objects are elements of $P$; there is a single morphism $p\to q$ if and only if $p\leqslant q$, otherwise there is no morphism.
            That is,
            \begin{equation*}
                \Hom(p,q)=
                \begin{cases}
                    \{p\leqslant q\} &\mbox{if } p\leqslant q;\\
                    \varnothing &\mbox{if } p>q.
                \end{cases}\qedhere
            \end{equation*}
        \end{definition}
        By definition all poset categories are small.
        One particular example that we often use is $\rrleq$, where $\leqslant$ is the usual partial ordering on $\rr$.

        \begin{table}[ht]
            \centering
            \begin{tabular}{lll}
                Name & Objects & Morphisms\\
                \toprule
                $\mathsf{Vec}$ & finite-dimensional vector spaces over $\zz_2$ & linear maps\\
                $\mathsf{Vec}_\infty$ & all vector spaces over $\zz_2$ & linear maps\\
                $\Top$ & topological spaces & continuous maps\qedhere
            \end{tabular}
            \caption{Definitions of some relevant categories.}\label{tb:useful-categories}
        \end{table}

        Our main interest is in $\Vecrrdiag$ diagrams\footnote{
            Though many that we come across can be factored through $\Top$, i.e. can be written as a composition of a $[\rrleq,\Top]$ diagram and a $[\Top,\mathsf{Vec}]$ diagram
        }, of which there are two relatively well-behaved classes: \emph{tame} and \emph{finite type} (though it later turns out that these are actually equivalent -- see \cref{le:tame-iff-finite-type}).

        \begin{definition}[Characteristic diagram]
            Let $I\subseteq\rr$ be an interval.
            Define the \emph{characteristic diagram $\chi_I\in\Vecrrdiag$} by
            \begin{align*}
                \chi_I(a)&=
                \begin{cases}
                    \zz/2\zz &\mbox{if }a\in I;\\
                    0 &\mbox{otherwise}.
                \end{cases}\\
                \chi_I(a\leqslant b)&=
                \begin{cases}
                    \id_{\zz/2\zz} &\mbox{if }a\in I;\\
                    0 &\mbox{otherwise}.
                \end{cases}\qedhere
            \end{align*}
        \end{definition}

        These characteristic diagrams behave nicely with finite intervals, and we can simplify things quite easily.
        Two useful examples are as follows.
        \begin{align}
            \chi_{[a,b)}\oplus\chi_{[b,c)}=\chi_{[a,c)}\quad&\text{for }a<b<c\in\rr;\label{eq:characteristic-diagram-reduction-1}\\
            \chi_I\oplus\chi_J=\chi_{I\cup J}\oplus\chi_{I\cap J}\quad&\text{for }I\cap J\neq\varnothing.\label{eq:characteristic-diagram-reduction-2}
        \end{align}

        \begin{definition}[Critical values]
            Let $F\in\Vecrrdiag$, and $I\subseteq\rr$ be an interval.
            We say that $F$ is \emph{constant on $I$} if $F(a\leqslant b)$ is an isomorphism for all $a\leqslant b\in I$.
            We say that $a\in\rr$ is a \emph{regular value} of $F$ is there exists some \emph{open} interval $J\subseteq\rr$ with $a\in J$ such that $F$ is constant \mbox{on $J$}.
            If $a\in\rr$ is not a regular value then we say that it is a \emph{critical value}.
        \end{definition}



        \begin{definition}[Finite type and tameness]
            Let $F\in\Vecrrdiag$.
            We say that $F$ is of \emph{finite type} if there exist finitely many intervals $I_1,\ldots,I_N\subseteq\rr$ such that $F=\bigoplus_{i=1}^N\chi_{I_i}$ and that $F$ is \emph{tame} if it has finitely-many critical values.
        \end{definition}

        \begin{lemma}\label{le:tame-iff-finite-type}
            Let $F\in\Vecrrdiag$.
            Then $F$ is tame if and only if $F$ is of finite type.
        \end{lemma}

        \begin{proof}
            \cite[Theorem~4.6]{Bubenik:dn}
        \end{proof}

        It turns out that we can define a notion of distance (though not quite a metric) between $\rrleq$-index diagrams, which provides useful when we start looking at applications of $\Vecrrdiag$ to algebraic topology.
        First, though, we need some more technical machinery.

        \begin{definition}[Translation functors and translation natural transformations]
            Let $b\geqslant0$.
            Define the \emph{$b$-translation functor $T_b$} by
            \begin{align*}
                T_b\colon\rrleq&\to\rrleq\\
                a&\mapsto a+b
            \end{align*}
            and define the \emph{$b$-translation natural transformation $\eta_b\colon\id_{\rrleq}\nt T_b$} by
            \begin{equation*}
                \eta_b(a) = a\leqslant a+b.\qedhere
            \end{equation*}
        \end{definition}

        It follows straight from the definitions that $T_bT_c=T_{b+c}$ and $\eta_b\eta_c=\eta_{b+c}$.

        \begin{definition}[Interleaving of diagrams]
            Let $F,G\in[\rrleq,\dcat]$ for some arbitrary category $\dcat$, and let $\vare\geqslant0$.
            Define an \emph{$\vare$-interleaving of $F$ and $G$} as a quadruple $(F,G,\varphi,\psi)$, where $\varphi\colon F\nt GT_\vare$ and $\psi\colon G\nt FT_\vare$ are natural transformations such that
            \begin{align*}
                (\psi T_\vare)\varphi &= F\eta_{2\vare},\\
                (\varphi T_\vare)\psi &= G\eta_{2\vare}.
            \end{align*}
            That is, we want the following diagrams to commute:
            \begin{equation*}
                \begin{tikzcd}[column sep=small]
                    F(a)\arrow{rr}{\eta_{2\vare}(a)}\arrow{dr}[swap]{\varphi(a)} & & F(a+2\vare) & & F(a+\vare)\arrow{dr}{\varphi(a+\vare)} &\\
                    & G(a+\vare)\arrow{ur}[swap]{\psi(a+\vare)} & & G(a)\arrow{rr}[swap]{\eta_{2\vare}(a)}\arrow{ur}[swap]{\psi(a)} & & G(a+2\vare)
                \end{tikzcd}
            \end{equation*}
            We say that $F$ and $G$ are \emph{$\vare$-interleaved} if there exists some $\vare$-interleaving $(F,G,\varphi,\psi)$.
        \end{definition}

        \begin{definition}[Interleaving extended pseudometric]\label{df:interleaving-metric}
            Let $\dcat$ be some arbitrary category.
            Define the extended pseudometric\footnote{
                See \cite[Theorem~3.3]{Bubenik:dn}; we quote: \emph{``[i]t fails to be a metric because it can take the value $\infty$ and $\dist(F,G)$ does not imply that $F\cong G$''}.
            } $\dist$ on any subset of the class of $[\rrleq,\dcat]$ diagrams by
                \begin{equation*}
                    \dist(F,G)=\inf_{\vare\geqslant0}\{\vare \mid F\text{ and }G\text{ are }\vare\text{-interleaved} \}.\qedhere
                \end{equation*}
        \end{definition}

        We now define one of the fundamental concepts in computational algebraic topology: \emph{persistent homology}.
        A good introduction to how this seemingly abstract definition arises in a reasonably natural way can be found in \cite[\S\S~5.13~\&~7.2]{ghrist2014elementary}, and our definition is from \cite[\S2.2.4]{Bubenik:dn}.

        \begin{definition}[Persistent homology]\label{df:persistent-homology}
            Let $F\in[\rrleq,\Top]$.
            Define the \emph{$p$-persistent $k$-th homology group $P_pH_kF(a)$ of $F$ at $a$} to be the image of the homomorphism $H_kF(a\leqslant a+p)$.
        \end{definition}

        This definition of persistent homology is better explained after some unpacking.
        Let $a,p\in\rr$, $k\in\nn$, and $F\in[\rrleq,\Top]$.
        Write $X_b$ to mean $F(b)$.
        Then
        \begin{itemize}
            \item $F(a\leqslant a+p)\colon X_a\hookrightarrow X_{a+p}$ is an inclusion map of topological spaces, since $F$ is a functor $\rrleq\to\Top$;
            \item $H_kF(a\leqslant a+p)\colon H_kX_a\to H_kX_{a+p}$ is the induced homomorphism of homology groups (which are $\zz_2$-vector spaces);
            \item $P_pH_kX_a=\im H_kF(a\leqslant a+p)\leqslant H_kX_{a+p}$ is a subgroup (subspace) of the $k$-th homology group of $X_{a+p}$.
        \end{itemize}
        Putting this all together\footnote{
            And using the fact that any topological space $X$ can be written in the form $F(a)$ for some $F\in[\rrleq,\Top]$ and $a\in\rr$ by taking the trivial diagram $F(a)=X$ for all $a\in\rr$.
        } we see that $P_pH_k\in[\Top,\mathsf{Vec}_\infty]$, where the functoriality follows by definition.
        Some more intuition behind the idea of persistent homology is given in \cref{sec:extended_persistence_homology}, where we explain the following statement:
        \begin{quotation}
            The $p$-persistent $k$-th homology group of $F$ at $a$ consists of homology classes that were \emph{born no later than $a$} and that are \emph{still alive at $a+p$}.
        \end{quotation}

        \subsection{Stability for persistent homology} 
        \label{sub:stability_for_persistent_homology}

            With these definitions and lemmas in hand, let us consider the following scenario: take some topological space $X$ and some (not necessarily continuous\footnote{
                See \cref{nt:why-not-assumed-continuous}.
            }) function $f\colon X\to(-M,M)\subset\rr$.
            We can define a \emph{height filtration $F\in[\rrleq,\Top]$ of $X$} by
            \begin{equation*}
                F(a)=f^{-1}\big((-\infty,a]\big)
            \end{equation*}
            and where $F(a\leqslant b)$ is the inclusion $F(a)\hookrightarrow F(b)$.
            Then we can define the $[\rrleq,\mathsf{Vec}_\infty]$ diagram $H_kF$ given by taking the $k$-th homology group of $F(a)$.
            For simplicity \emph{we assume}\footnote{
                See \cref{nt:compactness}.
            } \emph{that $H_kF$ is tame for all $k\in\nn$}.
            In particular then, by \cref{le:tame-iff-finite-type},
            \begin{equation}
                H_kF=\bigoplus_{i=1}^N\chi_{I_i},
            \end{equation}
            and so, for all $a\in\rr$,
            \begin{equation}
                H_kF(a)=\bigoplus_{i=1}^{N_a}(\zz/2\zz)
            \end{equation}
            for some $N_a\leqslant N$.
            That is, $H_kF$ being tame implies that all the homology groups $H_kF(a)$ are finitely generated, i.e. $H_kF\in\Vecrrdiag$.

            \begin{definition}[$M$-bounded tame functions]
                We call any such\footnote{
                    i.e. $f\colon X\to(-M,M)$ with $F(a)=f^{-1}\big((-\infty,a]\big)$ being such that $H_kF$ is tame.
                } $f$ an \emph{$M$-bounded tame function on $X$}, and call $F$ the \emph{associated filtration}\footnote{
                    Again, this is not standard terminology.
                }.
            \end{definition}

            \begin{note}[The tameness assumption]\label{nt:compactness}
                If we didn't assume that $H_kF$ is necessarily tame, but instead that $F(a)$ is a compact manifold for all $a\in\rr$, then all the (singular) homology groups $H_kF(a)$ are still finitely generated\footnote{
                    As a statement in full generality, this is reasonably non-trivial: see \cite[Proposition~III.1,~p.~130]{Neeb:2004te}
                }.
                So assuming that $H_kF$ is tame is at least no stronger than assuming that all of our topological spaces are compact manifolds -- this gives us a vague lower bound for the level of generality at which we are working.
            \end{note}

            Now, as in \cref{df:persistent-homology}, we can look at $p$-persistent homology.
            Since each $H_kF(a+p)$ is finitely generated, and we are working with $\zz_2$ coefficients\footnote{
                Here the fact that homology groups are actually vector spaces is vital, since it is a simple fact that the subspace of a finite-dimensional vector space is itself finite dimensional.
                If, however, we were working with coefficients in a general group then we would have to appeal to something like Schreier's lemma, which tells us that any \emph{finite index} subgroup of a finitely-generated group is itself finitely generated, or maybe even to some similar property of modules over a PID.
            }, the subgroup $P_pH_kF(a)$ is also finitely generated.
            That is, $P_pH_k\in[\Top,\mathsf{Vec}]$ and so $P_pH_kF\in\Vecrrdiag$.
            This means that we can use the interleaving extended pseudometric from \cref{df:interleaving-metric} on $p$-persistent homology groups of $F(a)$.

            \bigskip

            One of the main examples of an $M$-bounded tame function $f$ on a topological space $X$ is a \emph{height function}: we immerse $X$ into $\rr^n$ for some $n$ and `measure' $X$ along some axis\footnote{
                This is exactly the sort of example that we look at in \cref{sub:the_klein_bottle_persistent}; we explain how height functions relate to the assumption that $f$ is not necessarily continuous in \cref{nt:why-not-assumed-continuous}.
            }.
            In this case, we can obtain different height functions, and thus different associated $[\rrleq,\Top]$ diagrams, simply by perturbing the axis along which we measure by some small amount.
            For persistent homology to have much practical use we would strongly desire that small perturbations of the height functions result in small changes to the persistent homology groups.
            Explicitly, we would hope to be able to bound the distance\footnote{
                Measured by the interleaving extended pseudometric $\dist$.
            } between the persistent homology groups of $F$ and $G$ by the distance\footnote{
                The most natural choice of metric for data sampling being the sup metric $\|f-g\|_\infty$.
            } between $f$ and $g$.
            It turns out that this is, in fact, possible.

            \begin{theorem}[Stability theorem for persistent homology]\label{th:stability-for-persistent}
                Let $f,g$ be $M$-bounded tame functions on some topological space $X$, with associated filtrations $F,G$ (respectively).
                Then\footnote{
                    Recall that $\|f-g\|_\infty=\sup_{x\in X}|f(x)-g(x)|$.
                }
                \begin{equation*}
                    \dist(P_pH_kF,P_pH_kG)\leqslant\|f-g\|_\infty.\qedhere
                \end{equation*}
            \end{theorem}

            \begin{proof}
                \emph{By our previous comments -- namely that $P_pH_k\in[\Top,\mathsf{Vec}]$ -- this is a specific case of \cite[Theorem~5.1]{Bubenik:dn}.}
                \emph{As such, a full proof can be found there; we give here a short sketch of the proof.}

                Let $\vare=\|f-g\|_\infty$.
                Then
                \begin{equation*}
                    F(a) = f^{-1}\big((-\infty,a]\big) \subseteq g^{-1}\big((-\infty,a+\vare]\big) = G(a+\vare),
                \end{equation*}
                and similarly $G(a)\subseteq F(a+\vare)$.
                Combining these gives us inclusions
                \begin{equation*}
                    F(a)\hookrightarrow G(a+\vare)\hookrightarrow F(a+2\vare)
                \end{equation*}
                which is, by definition, the same as the inclusion $F(a)\hookrightarrow F(a+2\vare)$.
                Similarly we have $G(a)\hookrightarrow F(a+\vare)\hookrightarrow G(a+2\vare)$.
                Thus $F$ and $G$ are $\vare$-interleaved.
                But then the functoriality of $P_pH_k$ ensures that $P_pH_kF$ and $P_pH_kG$ are $\vare$-interleaved (see \cite[Proposition~3.6]{Bubenik:dn}) which gives the required result.
            \end{proof}

            \Cref{th:stability-for-persistent} tells us that small perturbations to our `measuring' function result in small perturbations to the resulting $p$-persistent homology groups.
            However, if we are given some topological space $X$, or construct one from a data point cloud, the stability theorem does \emph{not} ensure that picking \emph{any} function $f\colon X\to\rr$ will result in persistent homology necessarily telling us anything useful about the space.

            Many different applications of the stability theorem can be found in \cite[\S6]{Edelsbrunner:2008gf}.
            Two clearly important example (explained in full detail in \cite[\S4]{CohenSteiner:2007is}) that stand out, however, are that of \emph{homology inference}: computing the homology of a space bound by a smooth surface by computing the homology arising from a finite sample of points from the space; and \emph{shape comparison}: using persistent homology to measure how similar two topological spaces embedded in $\rr^n $are.

            The key point behind both of these examples is that, although \cref{th:stability-for-persistent} is phrased in terms of two functions on the same topological space, we can actually use it for analysing the persistent homology of two \emph{different} spaces: given some $X$ embedded in $\rr^n$ we can define $\dist_X\colon\rr^n\to\rr$ by $\dist_X(z)=\inf_{x\in X}\|z-x\|_{\rr^n}$.
            If we have another space $Y$ with $\dist_Y$ defined similarly then we can apply the stability theorem\footnote{
                After restricting to some compact subset of $\rr^n$ containing both $X$ and $Y$, say.
            } to $\dist_X$ and $\dist_Y$ to bound the `homological differences' between $X$ and $Y$ by the `Euclidean-distance differences'\footnote{
                It turns out that we can actually then bound the `Euclidean-distance difference' $\|\dist_X-\dist_Y\|_\infty$ by the \emph{Hausdorff distance} between $X$ and $Y$.
                See \cite[\S6]{Edelsbrunner:2008gf}.
            } between $X$ and $Y$.

            \begin{note}
                Our summaries of homology inference and shape comparison are very brief, and thus skip over some of the finer, but very important, details.
                The subtleties are explained fully in \cite[\S4]{CohenSteiner:2007is}, but the main problem is that similar barcodes don't necessarily imply similar spaces, and vice versa.
                To quote,
                \begin{quotation}
                    \emph{Perhaps unexpectedly, the homology groups of $X^{+\delta}$ can be different from those of $X$, even when $X$ has positive homological feature size and $\delta$ is arbitrarily small.}

                    \ldots
                    
                    \emph{In particular, two shapes whose persistence diagrams are close are not necessarily approximately congruent.}
                \end{quotation}
                (It does turns out, however, that the `pathological behaviour' behind the first part of the quote actually almost never occurs in practice, and using functions that aren't simply distance functions can solve the problem in the second part of the quote.)
            \end{note}

\begin{figure}[ht]
                \centering
                \includegraphics[width=.28\textwidth]{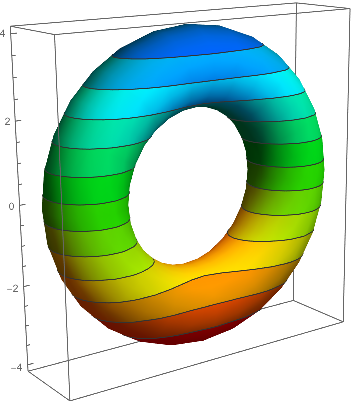}
                \includegraphics[width=.2\textwidth]{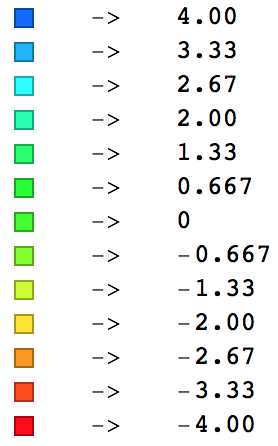}
                \hspace{.02\textwidth}
                \includegraphics[width=.28\textwidth]{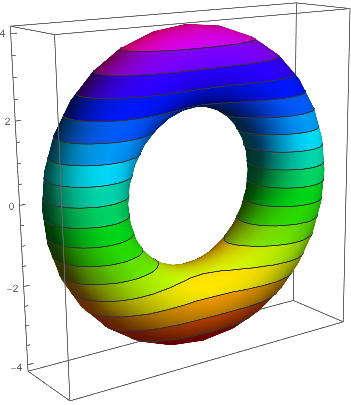}
                \includegraphics[width=.15\textwidth]{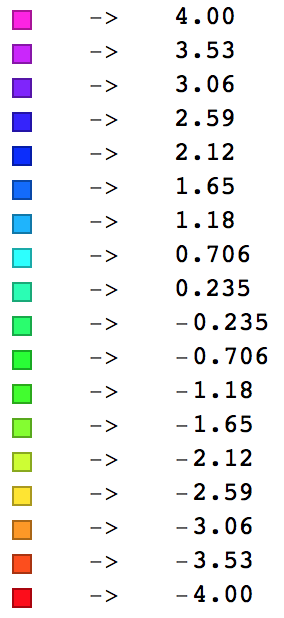}
                \caption{Discrete height functions on the torus using 13 and 18 partitions.}\label{fg:height-torus}
            \end{figure}

            \begin{note}[The continuity non-assumption]\label{nt:why-not-assumed-continuous}
                The fact that we don't require $f$ to be continuous corresponds to the idea that we might be using some sort of discrete height map, maybe because we are working with a simplicial complex or some other discrete version of our topological space.
                See \cref{fg:height-torus} for examples.
            \end{note}

            \bigskip

            If we have some $[\rrleq,\mathsf{Vec}]$ diagram of finite type then we can represent them graphically using \emph{barcodes}.
            Generally, barcodes give a very useful way of interpreting \cref{th:stability-for-persistent}, especially in light of \cite[Propositions~4.12,~4.13]{Bubenik:dn} which relates the interleaving distance of characteristic diagrams to the distances between the endpoints of their associated intervals.
            This is something that will prove useful when we examine the Klein bottle in \cref{sub:the_klein_bottle_persistent}.
            \begin{definition}[Barcodes]\label{df:barcode}
                Let $D\in\Vecrrdiag$ be of finite type, so that $D=\bigoplus_{i=1}^N\chi_{I_i}$ for some intervals $I\subseteq\rr$.
                Define the \emph{barcode $B_D$ of $D$} to be the multiset\footnote{
                    i.e. a set where each element occurs with a multiplicity: $\{x,x,y\}$ and $\{x,y\}$ define the same set, but different multisets.
                }
                \begin{equation*}
                    B_D = \{I_1,I_2,\ldots,I_N\}.\qedhere
                \end{equation*}
            \end{definition}

            A different visual representation of persistent homology is a \emph{persistence diagram}.
            These are introduced and discussed extensively in \cite{CohenSteiner:2007is}, as are the implications of the stability theorem.
            We choose, however, to use barcodes, and refer the interested reader to other sources for persistence diagrams.


        \subsection{Height of the Klein bottle (persistent homology)} 
        \label{sub:the_klein_bottle_persistent}

\begin{figure}[hpt]
                \centering
                \includegraphics[width=\textwidth]{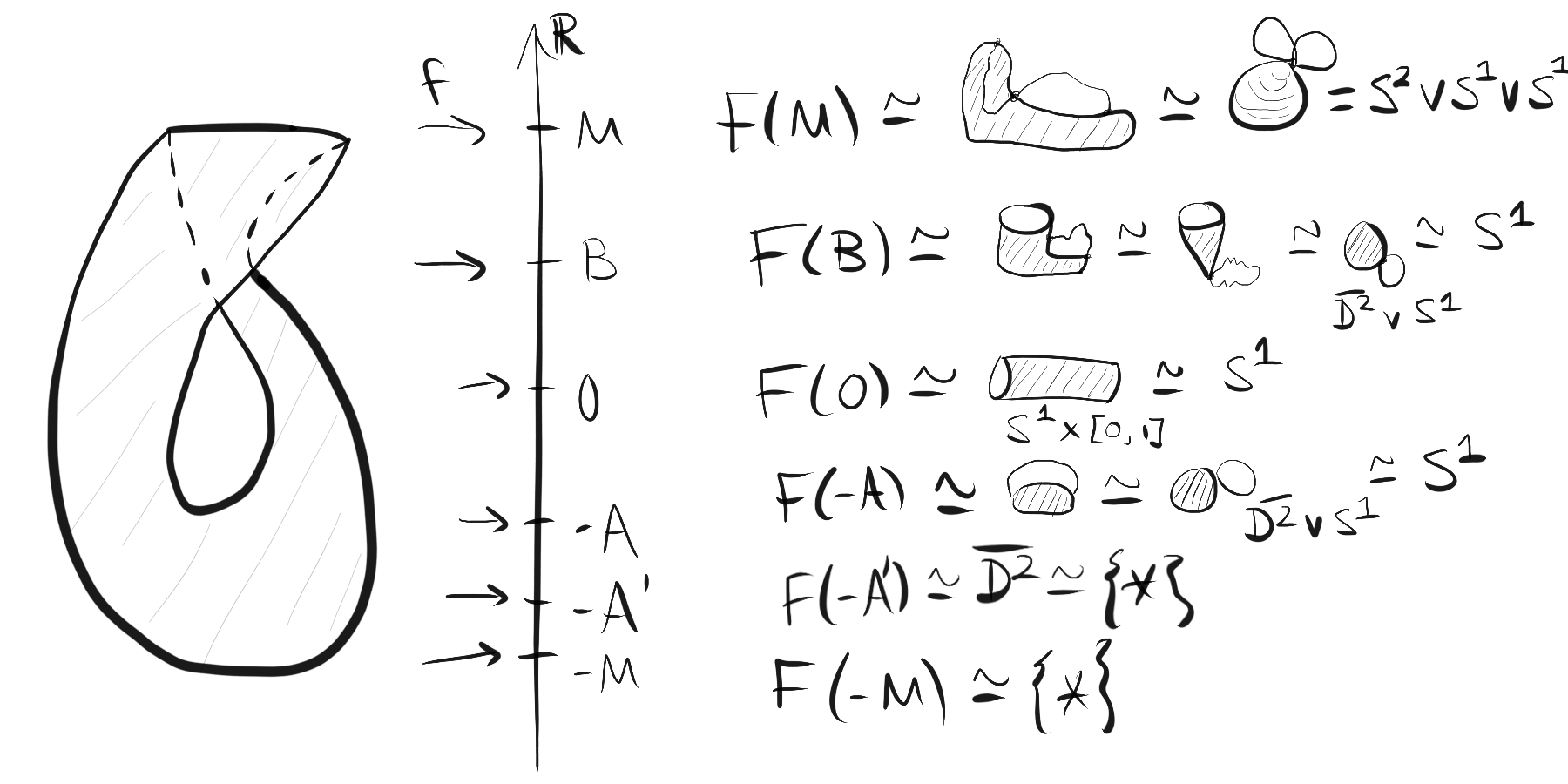}
                \caption{The Klein bottle $K$ immersed in $\rr^3$ with associated height function \mbox{$f\to(-M,M)\subset\rr$}. We shouldn't really be surprised that $F(-A)$, $F(0)$, and $F(B)$ are all homotopy equivalent, since the self-intersection at height $B$ isn't really an intersection at all: it is simply an artefact coming from this immersion in $\rr^3$.}\label{fg:klein-immersed-first-picture}
            \end{figure}

\begin{figure}[hpt]
                \centering
                \includegraphics[width=\textwidth]{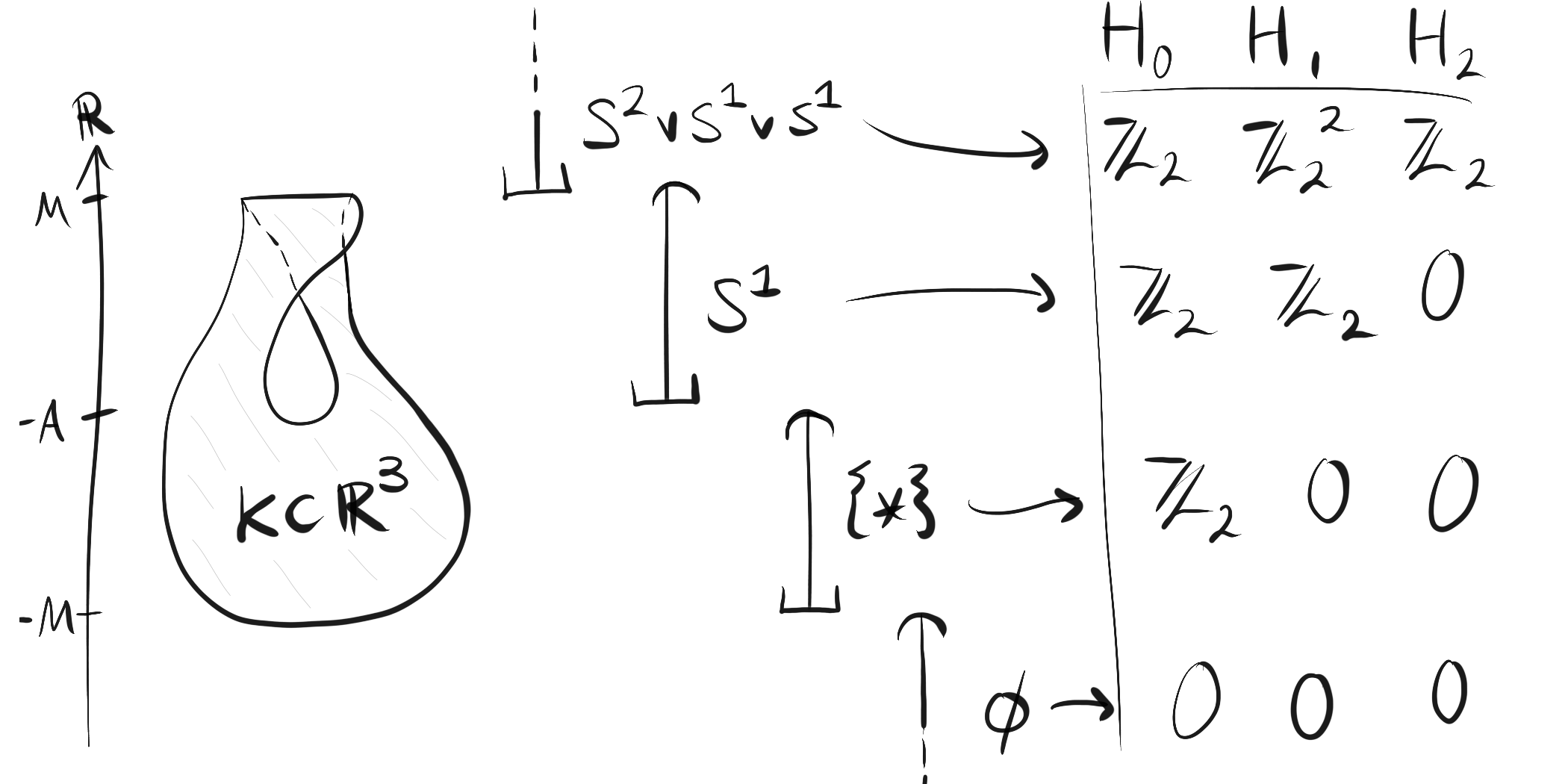}
                \caption{The homology groups of $F(a)$ between critical values.}\label{fg:klein-immersed-second-picture}
            \end{figure}

            We now study an explicit example: the Klein bottle, immersed in $\rr^3$.
            Here we scale the Klein bottle so that it has height $2M$ for some $M\in\rr$ and take $f$ to be the corresponding height function.
            Picking out certain values of $a\in\rr$ we can see how the `height slices' $F(a)$ change.
            We note that there are really three interesting points, namely $-M$, $A$, and $M$ (as labelled in \cref{fg:klein-immersed-first-picture}) where the homotopy type changes\footnote{
                We appeal to the fact that homotopy equivalent spaces have isomorphic (singular) homology groups (see \cite[Corollary~2.11,~\S2.1]{hatcher2002algebraic}).
            }, and so we can refine our picture\footnote{
                Where we calculate the homology groups using \cite[Corollary~2.25~\&~Proposition~A.5]{hatcher2002algebraic}.
            } (see \cref{fg:klein-immersed-second-picture}).
            Looking at the homology groups we can read the critical values of each $H_k$ straight off (see \cref{tb:critical-values}).

\begin{figure}[hpt]
                \centering
                \includegraphics[width=\textwidth]{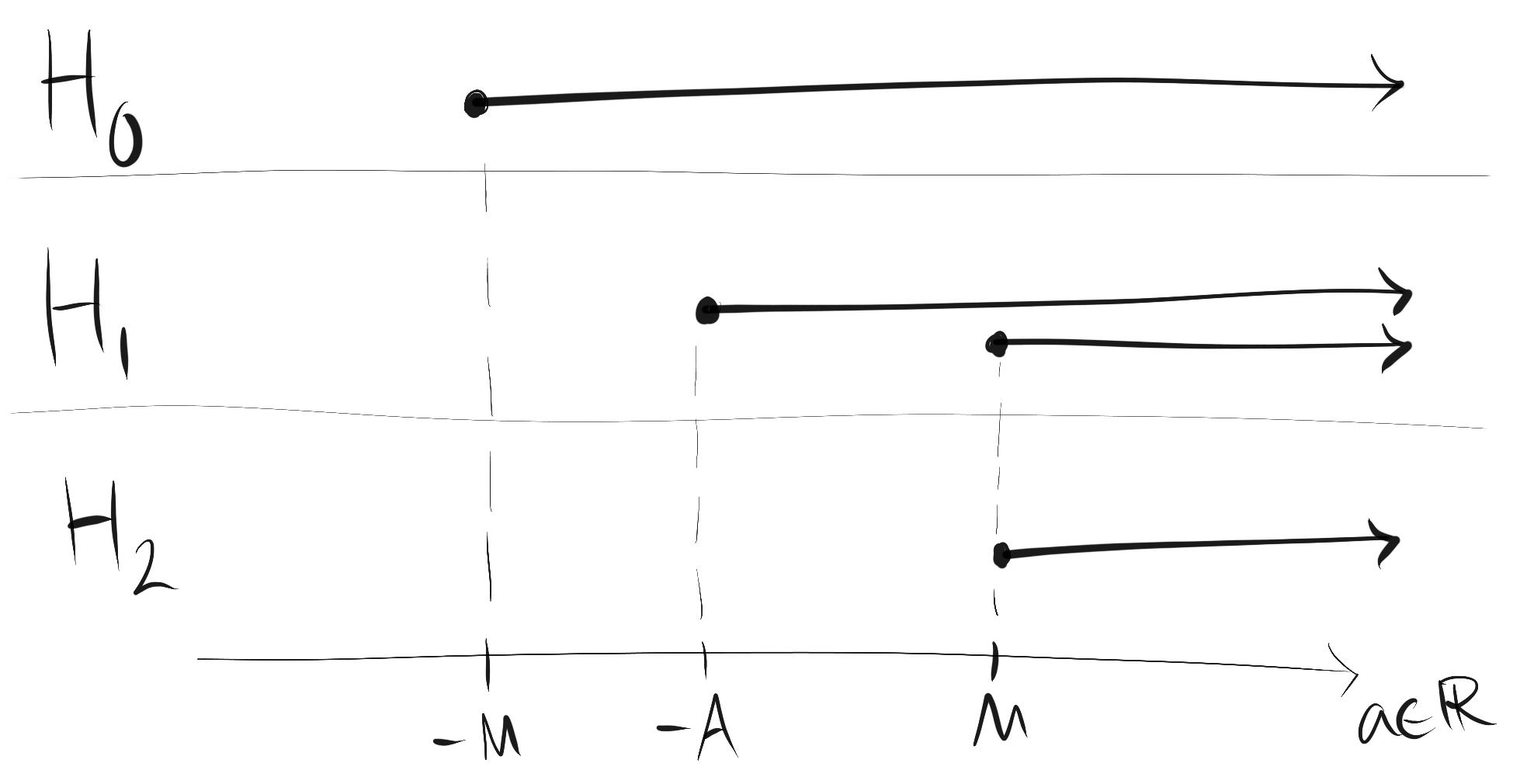}
                \caption{The barcode associated to the homology of $F(a)$. See \cref{eq:finite-type-decomopositions}.}\label{fg:klein-barcode}
            \end{figure}

            Putting all of the above together, we see that
            \begin{align}\label{eq:finite-type-decomopositions}
                H_0F &= \chi_{[-M,\infty)}\nonumber\\
                H_1F &= \chi_{[A,\infty)} \oplus \chi_{[M,\infty)},\\
                H_2F &= \chi_{[M,\infty)}.\nonumber
            \end{align}
            That is, each $H_KF$ is of finite type (and thus tame).
            We might have guessed this, since \cref{nt:compactness} told us that all of our homology groups would be finitely generated, but it is much simpler in this case to calculate finite-type decompositions for the $H_kF$ directly.
            We can summarise \cref{eq:finite-type-decomopositions} by using a barcode\footnote{
                Recall \cref{df:barcode}.
            } (see \cref{fg:klein-barcode}).
            Looking at the barcode we see one interesting feature: there are \emph{no deaths}\footnote{
                We explain the idea of \emph{birth} and \emph{death} in more detail in \cref{sec:extended_persistence_homology}.
            }.
            That is, as we increase $a$, at no point does the dimension of $H_kF(a)$ decrease.
            Alternatively, we can say that all of the intervals in our finite-type decomposition are upper-half infinite: if $a\in I$ then $a+p\in I$ for all $p\geqslant0$.
            So $H_kF(a\leqslant a+p)\cong\id_{H_kF(a)}$ by definition.
            That is, for all $p\geqslant0$ and $k\in\nn$,
            \begin{equation}
                P_pH_kF \cong H_kF.
            \end{equation}

\begin{figure}[hpt]
                \centering
                \includegraphics[width=\textwidth]{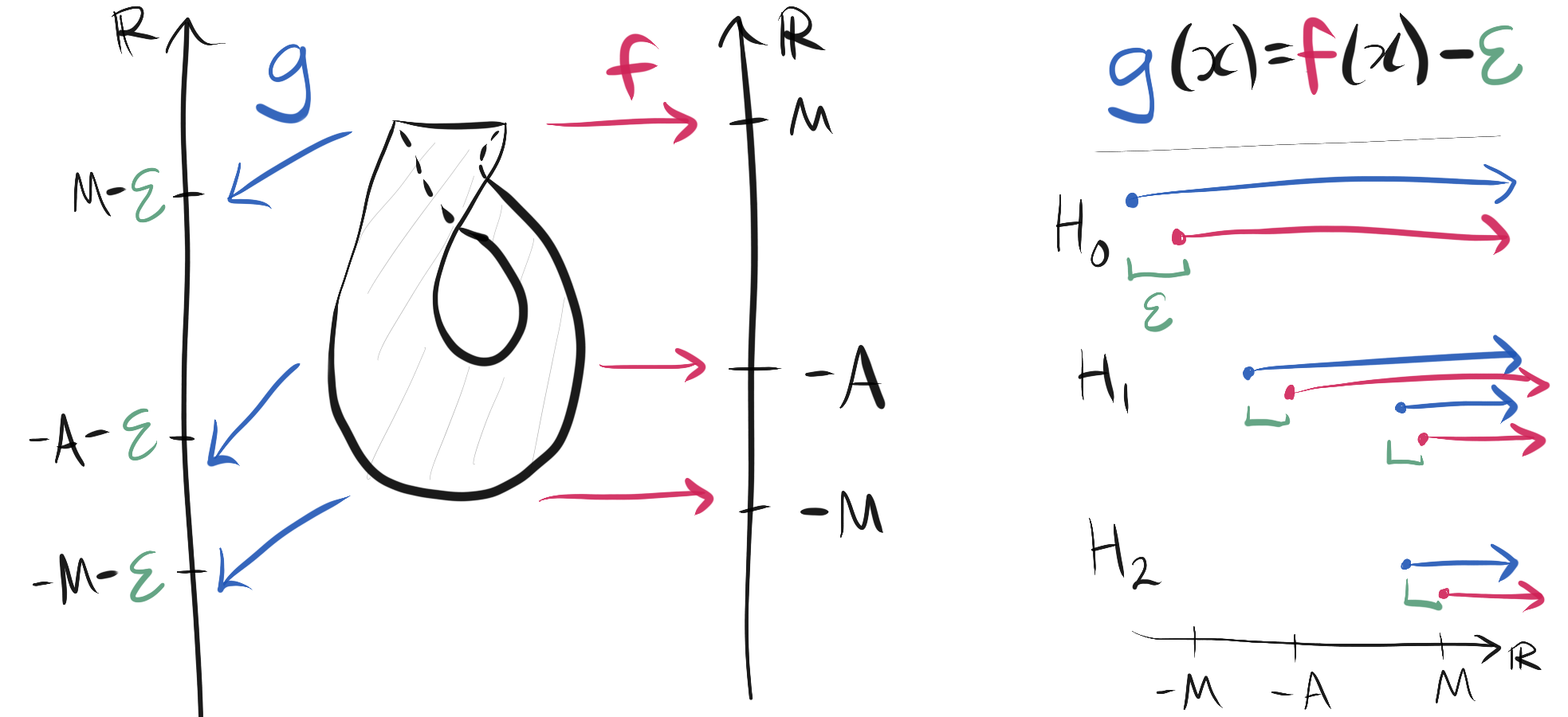}
                \caption{Here $g$ is simply $f$ shifted down by $\vare$. This means that $\|f-g\|_\infty=\vare$, and it turns out that $\dist(H_kF,H_kG)=\vare$ as well (we can precisely relate the interleaving distance to the \emph{bottleneck distance} between the barcodes, and there are some useful facts simplifying the calculation of interleaving distances between characteristic diagrams -- see \cite[Propositions~4.12,~4.13]{Bubenik:dn}).
                So the bound in \cref{th:stability-for-persistent} is attained (since here the persistent homology groups are exactly the homology groups -- we have no deaths).}\label{fg:shift-height-function-barcodes}
            \end{figure}

\begin{figure}[hpt]
                \centering
                \includegraphics[width=\textwidth]{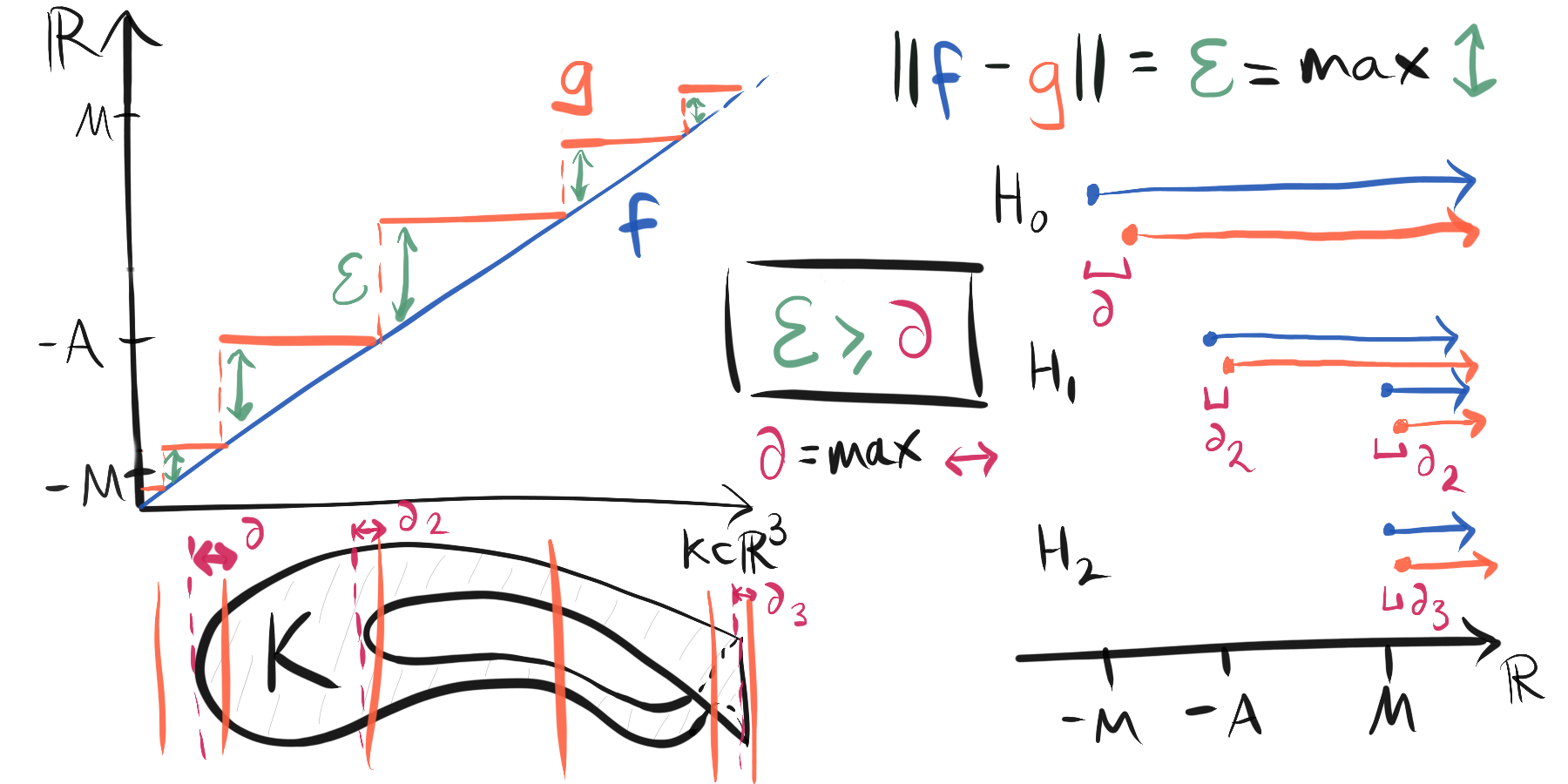}
                \caption{Using a discrete height function $g$ (as in \cref{fg:height-torus}). We see that, in this scenario, where the height of a section is determined by the highest value of $f$ taken on that section, the interleaving distance $\delta$ is given by the maximum distance between any partition boundary and the first critical point of $K$ occurring before it. This is bounded by $\|f-g\|_\infty$.}\label{fg:discrete-height-klein}
            \end{figure}

            Here then, \cref{th:stability-for-persistent}, along with \cite[Propositions~4.12,~4.13]{Bubenik:dn}, tells us that if $g$ is some other height function, close (in the $\|\cdot\|_\infty$ sense) to $f$, then the resulting barcodes\footnote{
                To be picky, we really mean the graphical representation of the barcodes, i.e. the associated intervals drawn inside $\rr$.
            } will be close (in the $|\cdot|$ sense).
            As a trivial example, we see that if we define $g$ as a shift of $f$ by $\vare$ then the resulting barcodes will be exactly distance $\vare$ apart (see \cref{fg:shift-height-function-barcodes}).
            As a slightly-less-contrived example, if we were working computationally and needed to use discrete data, then we could use a discrete height function (as in \cref{fg:height-torus}) and bound the errors on the resulting barcode (as compared to using a continuous height function) by how close the boundaries of our partitions are to critical points (see \cref{fg:discrete-height-klein}).
            Using the ideas of shape comparison, mentioned in \cref{th:stability-for-persistent}, we also know that deforming our immersion of the Klein bottle would result in bounded changes in the barcode.

            \begin{table}[hpt]
                \centering
                \begin{tabular}{lr}
                    $k$ & critical values of $H_kF$\\
                    \toprule
                    $0$ & $-M,M$\\
                    $1$ & $-A,M$\\
                    $2$ & $M$
                \end{tabular}
                \caption{Critical values of $H_kF$.}\label{tb:critical-values}
            \end{table}

    \section{Extended persistent homology} 
    \label{sec:extended_persistence_homology}

        \begin{quotation}
            \raggedleft
            \emph{``Death is the surest calculation that can be made.''}\\
            -- Ludwig Büchner, \emph{Force and Matter}
        \end{quotation}

        \emph{In this section we keep the same assumptions as in \cref{sec:persistent_homology}, namely that $f$ is some $M$-bounded tame function on a topological space $X$, and $F(a)=f^{-1}\big((-\infty,a]\big)$.}

        \bigskip

        Before moving on to discuss extended persistent homology, we first look about another way of talking at persistent homology.
        One very useful way of talking about persistent homology is using the idea of \emph{birth} and \emph{death}.
        This is really just a way of formalising some of the things that we noticed in \cref{sub:the_klein_bottle_persistent}.

        \bigskip

        For all $a\leqslant b$ we have the inclusion $\iota\colon F(a)\hookrightarrow F(b)$, which induces the homomorphism of homology groups $\iota_*\colon H_kF(a)\to H_kF(b)$.
        Using this, given some homology class $[\beta]\in H_kF(b)$ we can ask for which $a\leqslant b$ there exists some homology class $[\alpha]\in H_kF(a)$ with $\iota_*([\alpha])=[\beta]$.
        Clearly, if we find two such values of $a$, say $a_1$ and $a_2$, then it is the smaller one which is of the most interest: say $a_1\leqslant a_2$, then, using the general fact about induced homomorphisms that $(\iota\kappa)_*=\iota_*\kappa_*$, we can factor $H_kF(a_1)\to H_kF(b)$ through $H_kF(a_2)$.
        Thinking of $F(r)$ as evolving over time as $r\in\rr$ gets larger, we call the smallest such value\footnote{
            Which might well be $-\infty$.
        } of $a$ the \emph{time of birth of $[\beta]$}.
        Similarly, we can look at the \emph{time of death of $[\beta]$} by considering the image of $H_kF(a')\to H_kF(b')$ for $b'>b$ and all $a'<a$, where $a$ is the time of birth: it is the smallest\footnote{
            If any such $b'$ exists.
        } $b'$ such that $[\beta]$ is not in the image of $H_kF(a')\to H_kF(b')$.

        Using this language we can formulate the following motto of persistent homology:
        \begin{quotation}
            The $p$-persistent $k$-th homology group of $F$ at $a$ consists of homology classes that were \emph{born no later than $a$} and that are \emph{still alive at $a+p$}.
        \end{quotation}
        This makes it clear that, if we have no deaths, then persistent homology is simply homology -- every class that is born no later than $a$ will always be alive at $a+p$.

        \bigskip

        If $X$ is a surface then Morse theory tells us that the births and deaths of homology classes will be at critical points of the surface: if $[\beta]$ has birth time $a$ and death time $b$ then there will be critical points of $X$ at $F(a)$ and $F(b)$, call them $p_a$ and $p_b$, respectively.
        This gives us a pairing of critical points of $X$: we pair $p_a$ with $p_b$ and say that they have \emph{persistence} $|b-a|$ (or sometimes $|f(p_b)-f(p_a)|$).
        See \cite[\S2]{CohenSteiner:2009ho} for more details and motivation.

        The issue that remains (and that \cite{CohenSteiner:2009ho} aims to resolve) is that there are scenarios where homology classes don't die (i.e. have death time $\infty$), since this leaves some critical points of $X$ unpaired.
        \emph{Extended persistent homology} solves this problem by ensuring that every homology class eventually dies within finite time.

        \begin{definition}[Bifiltrations of an $M$-bounded tame function]\label{df:bifiltration}
            For an $M$-bounded tame function $f$ on a topological space $X$ define\footnote{
                The choice of the $+1$ in $\Fdown$ is arbitrary: we could use any `spacing' constant $\lambda>0$.
            } the \emph{associated $(2M+1)$-bifiltration}\footnote{
                Since we only work with $(2M+1)$-bifiltrations, we often refer to them just as \emph{bifiltrations}.
            } \emph{$\widehat{F}$} by
            \begin{align*}
                \Fup(a) = f^{-1}\big((-\infty,a]\big),\,&\,\,\,\Fdown(a) = f^{-1}\big([2M+1-a,\infty)\big),\\
                \widehat{F}(a)&=\Fpaira.
            \end{align*}
            So $\Fup=F$ in our previous notation, and where the $s$ in $\Fdown$ is to remind us that there is some shift, i.e. that $\Fdown(a)$ is not simply $f^{-1}\big([a,\infty)\big)$.
        \end{definition}

        The reason for these definitions is made slightly clearer when we look at how these functions\footnote{
            They are in fact diagrams, though this does require some reasoning, which we give later.
        } change as $a$ increases:
        \begin{equation}\label{eq:what-happens-to-fpair-as-a-grows}
            \widehat{F}(a)=
            \begin{cases}
                (\varnothing,\varnothing) &\mbox{for }a\in(-\infty,-M)\\
                (X_a,\varnothing) &\mbox{for }a\in[-M,M)\\
                (X,\varnothing) &\mbox{for }a\in [M,M+1)\\
                (X,X\setminus X_{2M+1-a}) &\mbox{for }a\in[M+1,3M+1)\\
                (X,X) &\mbox{for }a\in[3M+1,\infty)
            \end{cases}
        \end{equation}
        where $X_a=F(a)\subseteq X$ is a subspace of $X$.
        So we see that if we take the \emph{relative homology} of this pair $\widehat{F}(a)$ then we recover $H_kF(a)$ for \mbox{$a\in(-\infty,M+1)$}, since $H_k(Y,\varnothing)\cong H_k(Y)$ for all $Y$.
        But for $a\geqslant M+1$ the homology then `dies down', ending with all relative homology groups being $0$ for $a\geqslant 3M+1$, since $H_k(Y,Y)=0$ for all $Y$.
        See \cref{sub:the_klein_bottle_extended_persistent} for an example with the Klein bottle.

        The motivation for this construction of $\Fup$ and $\Fdown$ comes from \cite[\S6]{Bubenik:dn}, which is in turn motivated by the abstraction of the situation in \cite[\S4]{CohenSteiner:2009ho} where Poincaré and Lefschetz duality are used.
        We refer the reader to these two papers for further information; we carry on developing as much machinery as we can with the tools that we have.

        \bigskip

        We claim that $\Fup$ and $\Fdown$ are $\Vecrrdiag$ diagrams.
        That is, they are functorial: they preserve composition of morphisms and map identity morphisms to identity morphisms.
        This follows from \cref{eq:what-happens-to-fpair-as-a-grows}, since $\widehat{F}(a)\subseteq\widehat{F}(b)$ for all $a\leqslant b$, and this inclusion induces, in a functorial way, a homomorphism $H_k\widehat{F}(a)\to H_k\widehat{F}(b)$ on the relative homology groups\footnote{
            This works for relative homology almost exactly as it does for absolute homology, but there are some helpful comments just after Example~2.18 in \cite[118]{hatcher2002algebraic}.
        }.

        \begin{definition}[Extended persistent homology]
            Let $F\in[\rrleq,\Top]$.
            Define\footnote{
                Although we have been working with $M$-bounded tame functions $f$ on a topological space $X$, we can still define $\widehat{F}=(\Fup,\Fdown)$ for any $F\in[\rrleq,\Top]$ exactly as in \cref{df:bifiltration}.
            } the \emph{extended $p$-persistent $k$-th homology group $E_pH_k\widehat{F}(a)$ of $F$ at $a$} to be the image of the homomorphism $H_k\widehat{F}(a\leqslant a+p)$.
        \end{definition}

        As with \cref{df:persistent-homology}, this definition is better understood after some unpacking.
        Let $a,p\in\rr$, $k\in\nn$, and $F\in[\rrleq,\Top]$.
        Write $\widehat{X}_b$ to mean $\widehat{F}(b)$.
        Then
        \begin{itemize}
            \item $\widehat{F}(a\leqslant a+p)\colon\widehat{X}_a\hookrightarrow\widehat{X}_{a+p}$ is an inclusion map of topological spaces;
            \item $H_k\widehat{F}(a\leqslant a+p)\colon H_k\widehat{X}_a\to H_k\widehat{X}_{a+p}$ is the induced homomorphism of relative homology groups (which are $\zz_2$ vector spaces);
            \item $E_pH_k\widehat{X}_a=\im H_k\widehat{F}(a\leqslant a+p)\leqslant H_k\widehat{X}_{a+p}$ is a subgroup (subspace) of the $k$-th homology group of $\widehat{X}_{a+p}$.
        \end{itemize}
        Again, as with $P_pH_k$, we see that $E_pH_k\in[\Top,\mathsf{Vec}_\infty]$.

        We will see that, using \cref{eq:what-happens-to-fpair-as-a-grows}, we can sometimes think of extended persistent homology as follows.
        First we compute persistent homology `from bottom to top', then we compute persistent homology again, but from top to bottom and whilst squeezing our space to a point along the way.

        \subsection{Stability for extended persistent homology} 
        \label{sub:stability_for_extended_persistence_homology}

            \begin{theorem}[Stability theorem for extended persistent homology]\label{th:stability-for-extended-persistent}
                Let $f,g$ be $M$-bounded tame functions on some topological space $X$, with associated \mbox{$(2M+1)$-bifiltrations} $\widehat{F},\widehat{G}$, respectively.
                Then
                \begin{equation*}
                    \dist(E_pH_kF,E_pH_kG)\leqslant\|f-g\|_\infty.\qedhere
                \end{equation*}
            \end{theorem}

            \begin{proof}
                Let $\vare=\|f-g\|_\infty$.
                All we need to show is that $\widehat{F}$ and $\widehat{G}$ are $\vare$-interleaved, since then we can use \cite[Proposition~3.6]{Bubenik:dn}.
                As in the proof for \cref{th:stability-for-persistent}, this would follow from showing that $\widehat{F}(a)\subseteq\widehat{G}(a+\vare)$ and $\widehat{G}(a)\subseteq\widehat{F}(a+\vare)$ for all $a\in\rr$.

                If $a+2\vare\in(-\infty,M+1)$ then this follows as in the proof of \cref{th:stability-for-persistent}, since $\Fdown,G^\downarrow_s$ are both $\varnothing$ for $a$, $a+\vare$, and $a+2\vare$.
                Similarly, if $a\geqslant M+1$ then $\Fup(a)\subseteq G_\uparrow(a+\vare)\subseteq \Fup(a+2\vare)$, so all that remains to show is that $\Fdown(a)\subseteq G^\downarrow_s(a+\vare)\subseteq\Fdown(a+2\vare)$.
                But this follows from the observation that $\Fdown(a)=X\setminus\Fup(2M+1-a)$ (see \cref{eq:what-happens-to-fpair-as-a-grows}).
                For all other cases (say, $a+\vare\in(-\infty,M+1)$ but $a+2\vare\geqslant M+1$) we can combine the above two arguments to show the required inclusions.
            \end{proof}


        \subsection{Height of the Klein bottle (extended persistent homology)} 
        \label{sub:the_klein_bottle_extended_persistent}

\begin{figure}[hpt]
                \centering
                \includegraphics[width=\textwidth]{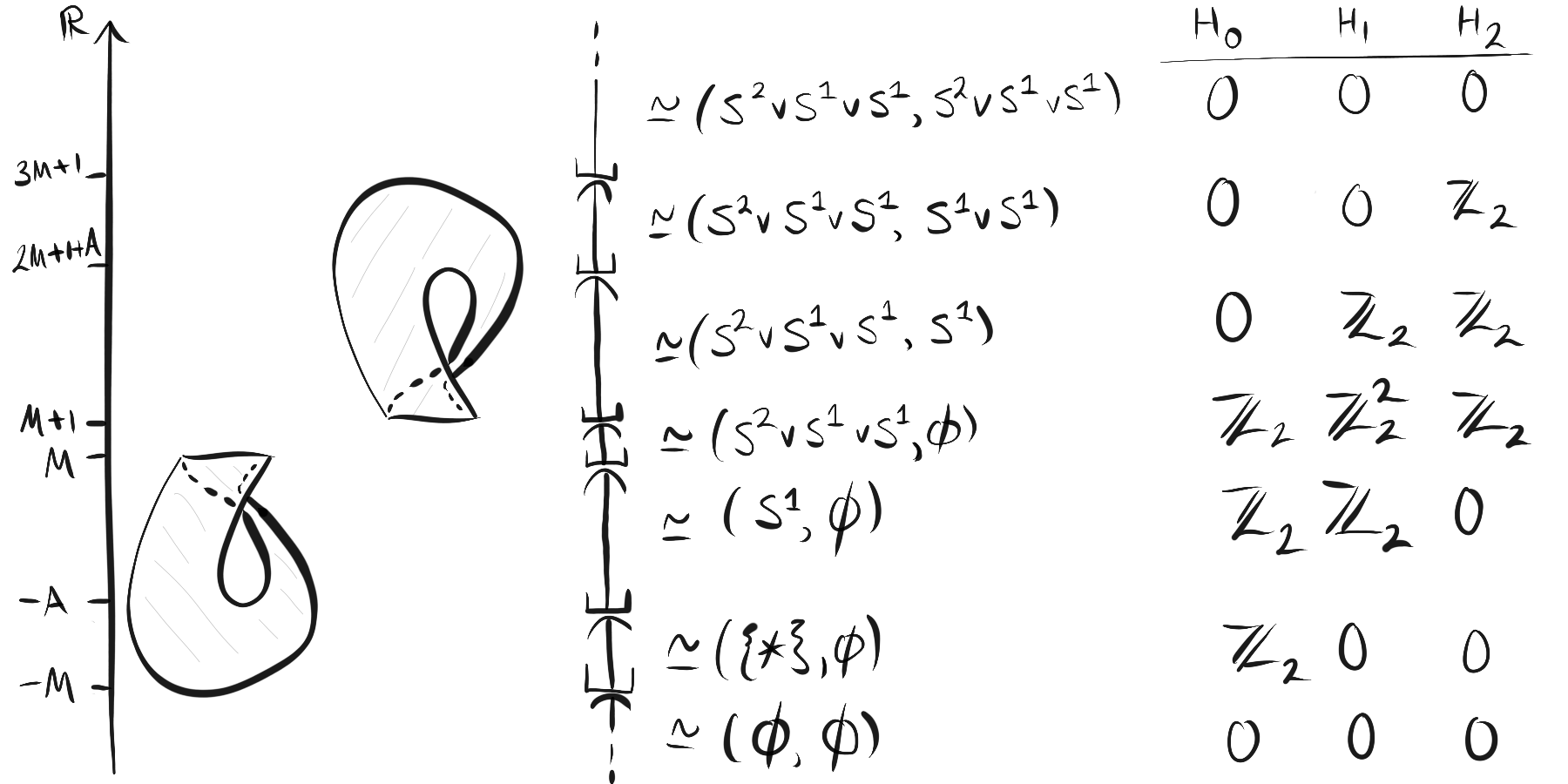}
                \caption{The homology groups and homotopy-equivalent spaces associated to the \mbox{$2M+1$}-bifiltration $\widehat{F}(a)$ of $K$. Note the reversed symmetry in the table of relative homology groups.}\label{fg:bifiltration-klein}
            \end{figure}

\begin{figure}[hpt]
                \centering
                \includegraphics[width=.8\textwidth]{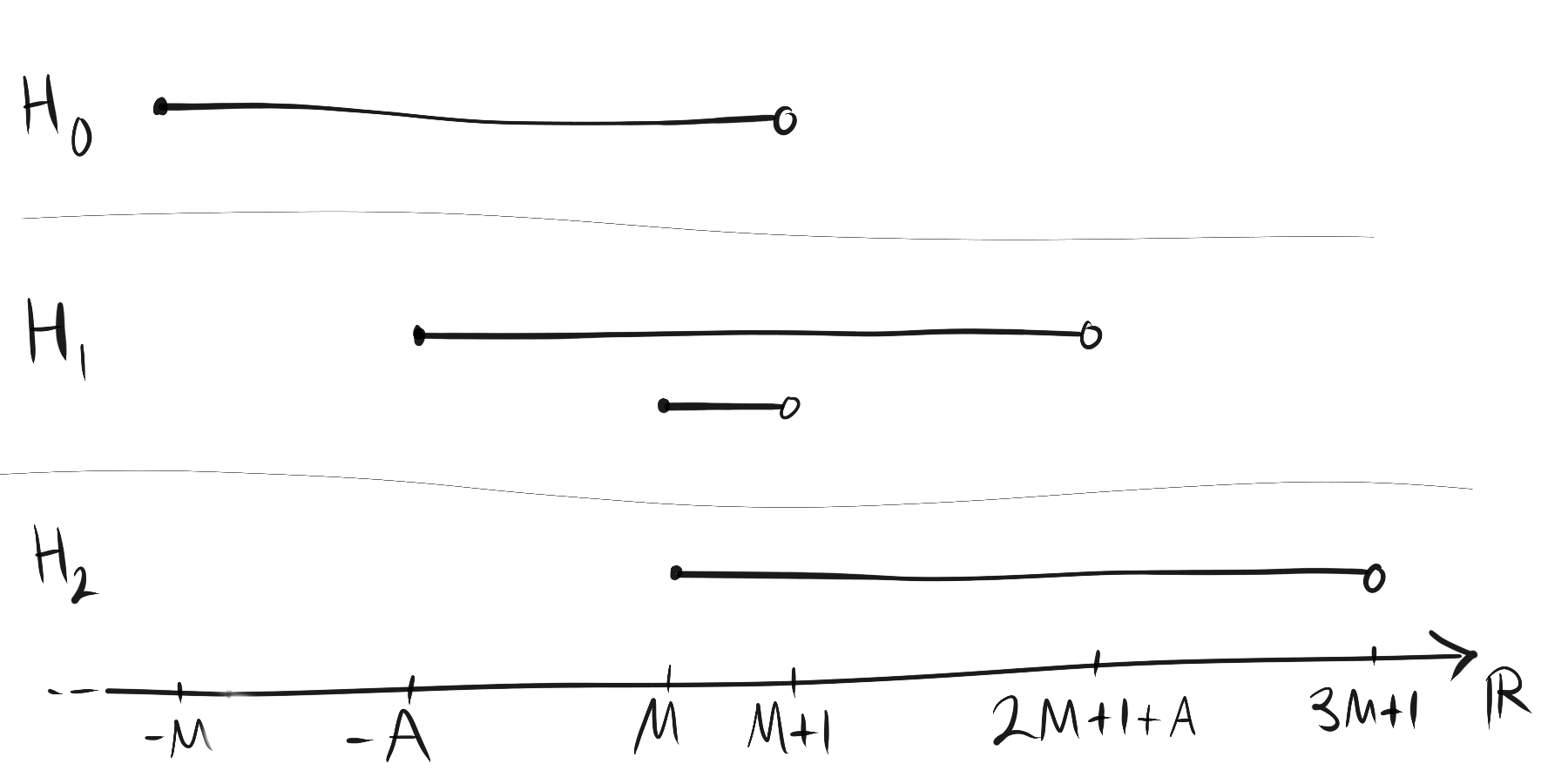}
                \caption{The barcode associated to the \mbox{$2M+1$}-bifiltration $\widehat{F}(a)$ of $K$. The symmetry in the table in \cref{fg:bifiltration-klein} is even more apparent here: $H_1$ is symmetric around the midpoint of $[M,M+1]$, and $H_0,H_2$ seem to mirror each other. Note that all intervals are left-closed and right-open (i.e. of the form $[x,y)$).}\label{fg:klein-bifiltration-barcode}
            \end{figure}

            We now return to the example of the Klein bottle immersed in $\rr^3$ of height $2M$ with height function $f$, as in \cref{sub:the_klein_bottle_persistent}.
            Previously, the fact that there were no deaths (i.e. that every homology class had infinite persistence) meant that persistent homology looked exactly like homology.
            But here, for the same example, extended persistent homology \emph{guarantees death in finite time}.
            This means that we will expect to see different results when we compute $E_pH_kF(a)$ as compared to $H_kF(a)$.
            In \cref{fg:bifiltration-klein} we sketch $\widehat{F}(a)$, recalling \cref{eq:what-happens-to-fpair-as-a-grows}, and calculate homotopy-equivalent spaces\footnote{
                As in absolute homology, homotopy equivalent \emph{topological pairs} have the same homology.
                In a sense, this is trivial if you adopt the Eilenberg-Steenrod axiomatic point of view, since this is one of the axioms.
                On the other hand, we do have to be slightly careful: if we take a space $X$ and two \emph{homeomorphic} (and thus homotopy equivalent) subspaces $A,B\subset X$ then it is not necessarily true that $H_k(X,A)\cong H_k(X,B)$.
                For this to hold we also need that the inclusions $A\hookrightarrow X$ and $A\to B\hookrightarrow X$ are homotopic, where $A\to B$ is a homotopy equivalence.
            } and their homology groups for various $\widehat{F}(a)$.
            To calculate the homology groups, we use two facts: $H_k(X,\varnothing)\cong H_k(X)$; and \cite[Proposition~2.22~\&~Proposition~A.5]{hatcher2002algebraic}, which says that if $(X,A)$ is a good pair then $H_k(X,A)\cong\widetilde{H}_n(X/A)$.
            From this, we can draw the barcode and read off a finite-type decomposition of the $H_k\widehat{F}(a)$ (see \cref{fg:klein-bifiltration-barcode}).

            If we now turn to the extended persistent homology groups we get some non-trivial results.
            Each group is parametrised by two variables: $a\in\rr$, which we think of as time, and $p\geqslant0$, which we think of as lifespan.
            For example, if $E_pH_k\widehat{F}(a)$ is non-zero then it means that there is some $k$-homology class alive at time $a$ that persists until at least time $a+p$.
            The extended $p$-persistent $k$-th homology groups of $\widehat{F}$ at $a$ are as follows.

            \begin{equation}\label{eq:eph-klein-0}
                E_pH_0\widehat{F}(a)=
                \begin{cases}
                    0 &\mbox{if }a\not\in [-M,~M+1);\\
                    0 &\mbox{if }p\geqslant(2M+1);\\
                    0 &\mbox{if }a\in[-M,~M+1)\text{ but }p\geqslant\big(2M+1-|a+M|\big);\\
                    \zz_2 &\mbox{otherwise}.
                \end{cases}
            \end{equation}

            \begin{equation}\label{eq:eph-klein-1}
                E_pH_1\widehat{F}(a)=
                \begin{cases}
                    0 &\mbox{if }a\not\in [-A,~2M+1+A);\\
                    0 &\mbox{if }p\geqslant \big(2(M+A)+1\big);\\
                    0 &\mbox{if }a\in[-A,~2M+1+A)\text{ but }p\geqslant\big(2(M+A)+1-|a+A|\big);\\
                    \zz_2\oplus\zz_2 &\mbox{if }a\in [M,M+1)\text{ and }p\leqslant(M+1-a);\\
                    \zz_2 &\mbox{otherwise}.
                \end{cases}
            \end{equation}

            \begin{equation}\label{eq:eph-klein-2}
                E_pH_2\widehat{F}(a)=
                \begin{cases}
                    0 &\mbox{if }a\not\in [M,~3M+1);\\
                    0 &\mbox{if }p\geqslant(2M+1);\\
                    0 &\mbox{if }a\in[M,~3M+1)\text{ but }p\geqslant(M+1-a);\\
                    \zz_2 &\mbox{otherwise}.
                \end{cases}
            \end{equation}

\begin{figure}[hpt]
                \centering
                \includegraphics[width=\textwidth]{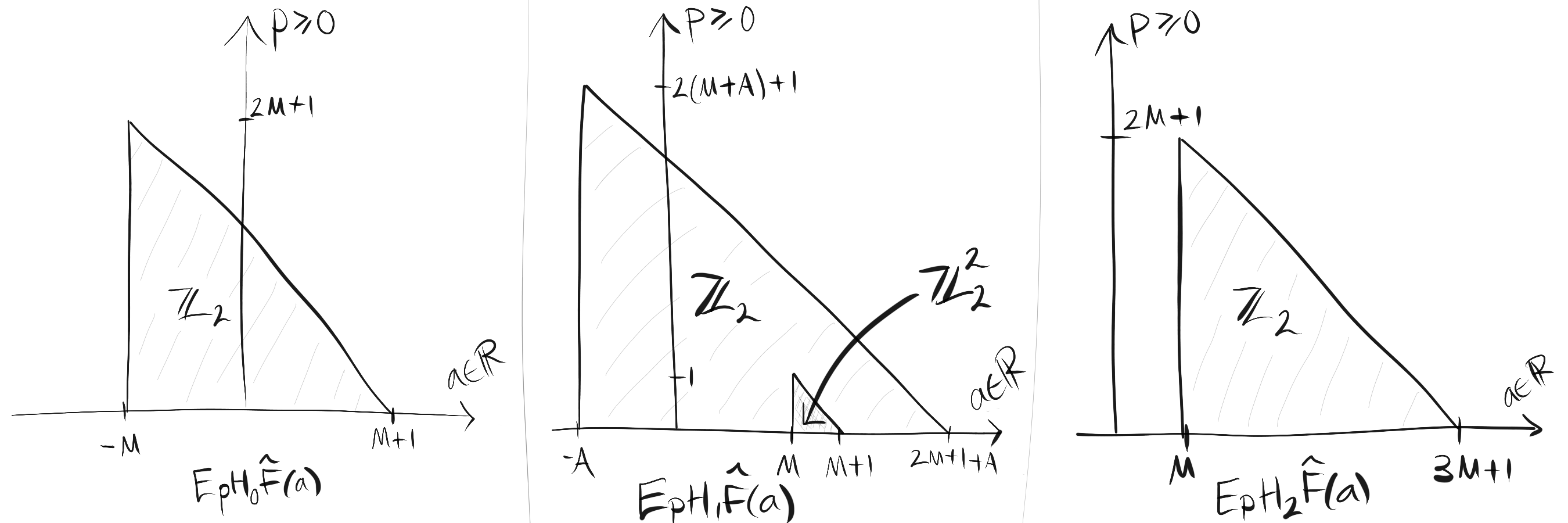}
                \caption{Extended persistent homology ensures death in finite time, and so each region in the $(a,p)$ plane of $E_pH_k\widehat{F}(a)$ is compact, as opposed to being half infinite.}\label{fg:2d-barcodes}
            \end{figure}

            
            Although these formulas might look slightly daunting at first, it is largely due to a lack of concise notation -- comparing them to \cref{fg:klein-bifiltration-barcode,fg:2d-barcodes} we see that they are (secretly) reasonably simple.
            Looking at the more general formula given in \cref{sub:extended_persistence_homology_for_general_tame_diagrams} will hopefully also help clarify what is actually being said in \cref{eq:eph-klein-0,eq:eph-klein-1,eq:eph-klein-2}.

            As in \cref{sub:the_klein_bottle_persistent}, we can use the stability theorem in many ways.
            One interesting point is that, since all of our intervals in the barcode are now finite, the stability theorem could now be interpreted as being slightly \emph{weaker} since extended persistence homology is naturally slightly \emph{stronger}: if we have two half-infinite intervals then their interleaving distance is exactly the distance between their two (finite) endpoints; if we have two \emph{finite} intervals then their interleaving distance is the \emph{minimum} of the (maximum) distance between their endpoints and their (maximum) length.
            That is, if we take a space, compute its persistent homology and extended persistent homology, both using two different functions $f$ and $g$, then the stability theorem tells us that the difference between the persistent homology barcodes will be no more than $\|f-g\|_\infty$.
            It tells us the same thing for the extended persistent homology barcodes, but we know that there is a distinct possibility that the difference between the barcodes will actually be smaller still.

            Of course, since this `double bound' consists of two `less-than-\emph{or-equal}' inequalities, it could very well be the case that there is no discernible improvement in this specific manner.
            However, it seems that it might be worth thinking about whether we can construct extended persistent homology in such a way that the barcode distance is \emph{strictly} less than that for persistent homology.

        \subsection{Extended persistent homology for general tame diagrams} 
        \label{sub:extended_persistence_homology_for_general_tame_diagrams}

            For the sake of completeness we now generalise \cref{eq:eph-klein-0,eq:eph-klein-1,eq:eph-klein-2} to a more general setting, though still assuming that our diagrams are tame.
            We start by looking at the simple case where $H_k\widehat{F}=\chi_I$, and then extend this to the full case where \mbox{$H_k\widehat{F}=\bigoplus_{i=1}^k\chi_{I_i}$} for some $k\in\nn$.

            \bigskip

            Clearly, if $H_k\widehat{F}=\chi_I$ for some interval $I\subset\rr$, then $E_pH_k\widehat{F}(a)=0$ whenever $a\not\in I$, since $H_k\widehat{F}(a)=0$.
            That is, any class that is not even alive has zero persistence.
            Similarly, if $p\geqslant|I|$ then $E_pH_k\widehat{F}(a)=0$ for any choice of $a$.
            That is, each homology group has some element that lives the longest, and so no class will persist longer than this maximum lifespan.

            In this simple case then, where $H_k\widehat{F}=\chi_I$, we can summarise the extended persistent homology groups quite neatly.
            Write $I=[s,t)$, where $t$ in necessarily finite.
            Then\footnote{
                This is not necessarily the simplest way of expressing $E_pH_k\widehat{F}(a)$ in terms of conditions on $a$ and $p$, but it is meant to give some intuition: the first condition corresponds to the fact that things that aren't born have zero persistence; the second to the fact that nothing can live longer than the maximal lifespan; and the third to the fact that, if something has already been alive for time $y$ then it must die after $T-y$ more time has passed, where $T$ is its total lifespan.
            }
            \begin{equation}\label{eq:general-one-interval-eph}
                E_pH_k\widehat{F}(a)=
                \begin{cases}
                    0 &\mbox{if }a\not\in [s,t);\\
                    0 &\mbox{if }p\geqslant|s-t|;\\
                    0 &\mbox{if }a\in[s,t)\text{ but }p\geqslant|s-t|-|a-s|;\\
                    \zz/2\zz &\mbox{otherwise}.
                \end{cases}
            \end{equation}

            In the more general (but still tame) case where $H_k\widehat{F}=\bigoplus_{i=1}^k\chi_{I_i}$ for intervals $I_i=[s_i,t_i)\subset\rr$, we can still obtain some general formula\footnote{
                Though, in practice, it is much easier to read this information straight off from the barcode.
            }.
            For $m\in\nn$, write $N_m=\{n\in\nn \mid n\leqslant m\}=\{1,2,\ldots,m\}$.
            Then\footnote{
                The last condition in \cref{eq:general-one-interval-eph} could be phrased in a different way: we say that the group is $\zz_2$, \emph{unless} $a\in I_{j_1},I_{j_2}$ for distinct $j_1,j_2$ \emph{and} $p$ is small enough that $a+p\in I_{j_1},I_{j_2}$, in which case the group is $\zz_2\oplus\zz_2$, \emph{unless} $a\in I_{j_1},I_{j_2},I_{j_3}$ for distinct $j_1,j_2,j_3$ \emph{and} $p$ is small enough that~\ldots.
                Obviously though, there are many different ways of phrasing a reasonably complicated set of if-then phrases; the phrasing in \cref{eq:general-one-interval-eph} was simply the first one that occurred to the author.
                In particular, it seems that a simpler form could be obtained by using \cref{eq:characteristic-diagram-reduction-1,eq:characteristic-diagram-reduction-2}.
            }
            \begin{align}
                E_pH_k\widehat{F}(a)=
                &\begin{cases}
                    0 &\mbox{if }a\not\in\bigcup_{i\in N_k}[s_i,t_i);\\
                    0 &\mbox{if }p\geqslant\max_{i\in N_k}|t_i-s_i|;\\
                    0 &\mbox{if }a\in I_j\text{ for some }j\in N_k\text{ but }p\geqslant|t_j-s_j|-|a-s_j|;\\
                    (\zz/2\zz)^d &\mbox{otherwise,}
                \end{cases}\nonumber\\
                &\qquad\quad\text{where }d=\min_{t\in(a,a+p)}\left\{|\sigma| : \sigma\subseteq N_k, t\in\bigcap_{i\in\sigma}I_i\right\}.
            \end{align}
            A point to note when calculating the above is that, even if at times $a$ and $a+p$ the homology group $H_k$ is non-zero, if it is zero at some time $t\in(a,a+p)$ then $H_k\widehat{F}(a\leqslant a+p)$ factors through zero, and so $E_pH_k\widehat{F}(a)$ is also zero.
            More generally, if $H_k$ is of dimension $d$ at times $a$ and $a+p$, if it drops dimension at some time $t\in(a,a+p)$ then $E_pH_k\widehat{F}(a)$ will be of dimension $d'=\min_{t\in(a,a+p)}\dim H_k\widehat{F}(t)$.

            However, in practice it is still much easier to simply draw the barcode and read the data straight off from there.
        


    \section{Practicality of computational homology inference} 
    \label{sec:persistent_homology_of_a_cow}

        \begin{quotation}
            \raggedleft
            \emph{``All the really good ideas I ever had\\ came to me while I was milking a cow.''}\\
            Grant Wood
        \end{quotation}

        We mentioned the idea of \emph{homology inference} from \cite[\S4]{CohenSteiner:2007is} previously in passing.
        Here we spend a small amount of time looking at the practicality of this method in terms of computation.

        \bigskip

        The stability theorem tells us that we can estimate the persistent homology of, for example, a smooth manifold, by looking at a finite discrete subset of points.
        We use the finite set of points to construct a \emph{Vietoris-Rips complex}: in essence, we consider a ball of radius $r$ around each point and introduce a simplex between vertices whose balls intersect.
        Letting $r$ vary we can obtain a (discrete) filtration, and so apply the techniques of persistent homology.
        
\begin{figure}[hpt]
            \centering
            \includegraphics[width=.32\textwidth]{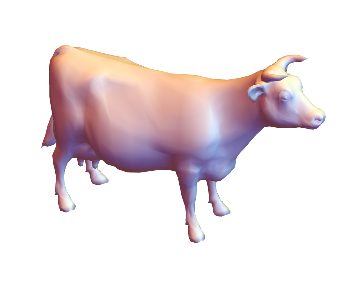}\includegraphics[width=.3\textwidth]{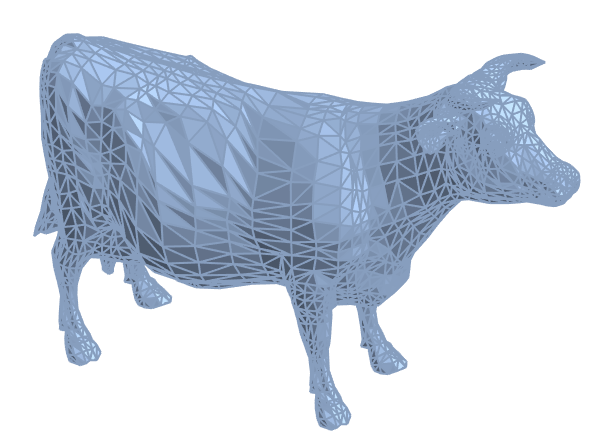}\\
            \includegraphics[width=.3\textwidth]{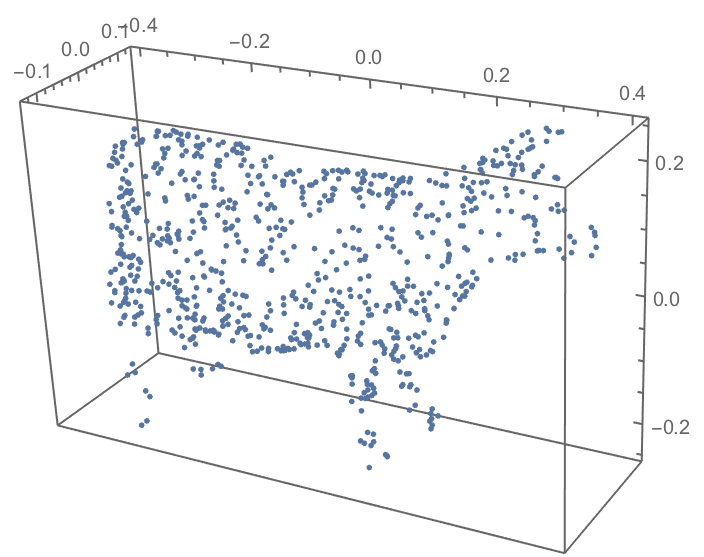}\includegraphics[width=.3\textwidth]{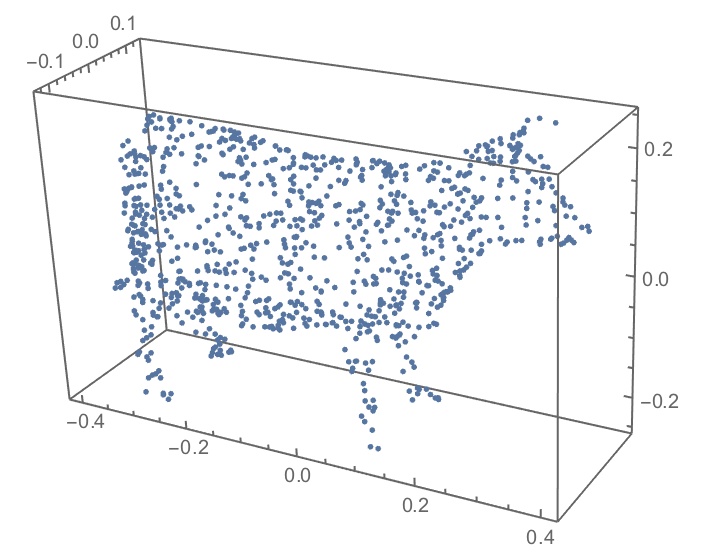}
            \caption{From top-left, clockwise: A 3D model of a cow; a discrete version of the model; a 1000-point random sample from the discrete model; an 800-point random sample from the discrete model.}\label{fg:cow}
        \end{figure}

        There is now quite a wide choice of software and libraries that can be used to compute persistent homology; \cite{Otter:2015wi} provides a thorough survey of the options available.
        Simply as a proof of concept though, we demonstrate here an example using \emph{Mathematica} \cite{Mathematica:tp} and \emph{Perseus} \cite{Perseus:70k_p2_v}.
        We pick a reasonably complex 3D model from the \texttt{ExampleData[``Geometry3D'']} library provided by Mathematica -- see \cref{fg:cow}.
        We generate a discrete version this model and use the \texttt{DirichletDistribution} function\footnote{
            The uniform distribution on any simplex.
        } to randomly sample points from the discrete model\footnote{
            Using a method given on \url{http://mathematica.stackexchange.com/questions/57938/} by user \texttt{ybeltukov}.
        }.
        We can then use the \texttt{brips} option of Perseus to calculate the persistent homology of the point cloud by using the Vietoris-Rips complex.
        There are two (relevant) options that we can alter: the amount $s$ by which the radius increases at each step, and the number $N$ of total steps to calculate.
        Using the Betti numbers (i.e. the dimensions of the homology groups, noting that Perseus also works with coefficients in $\zz/2\zz$) at each stage, we can import this data back into Mathematica and use \texttt{MatrixPlot} to generate a sort of barcode: instead of having multiple lines in each degree of homology, we have one line and represent the dimension by brightness -- the darker the segment the higher the dimension, with the key below showing explicit values.
        These barcodes are shown in \cref{tb:cow-barcodes}.

        The choice of taking an 800- and a 1000-point sample was arbitrary\footnote{
            It just so happens that 800 points was enough for the author to easily recognise the model, and 1000 was roughly the highest number of points for which Perseus could run in a few seconds with the given $(N,s)$.
        }, but we can still read some interesting data from the barcodes.
        For example, we see that, for $N=40,80$, the degree-$1$ homology dies off quicker in the 1000-point sample, and the degree-$2$ homology is born again at a later time (and so the persistent homology more closely resembles the homology of the original model).

        Although these barcodes aren't vastly different, that is to be expected, exactly because of the stability theorem: comparing the point samples in \cref{fg:cow} we see that the maximum distance between any two points in the 800 sample is not much more than the maximum distance between any two points in the 1000 sample, and so the difference in the barcodes will be small as well.
        In a sense, the stability theorem tells us that we can do homology inference with a finite set of points, but also that, if we choose our points in a uniform way so that they are roughly evenly distributed then any two random samples will give similar barcodes.
        An interesting consequence of this is that, if we know our point samples are roughly uniformly distributed, it might usually be enough to compute the persistent homology \emph{only once} -- running the computation again with different (uniformly distributed) data won't result in many new results.
        This is a point that we mention again, as well as the questions that it raises, in \cref{sec:conclusions}.

        \begin{table}[hpt]
            \centering
            \begin{tabular}{r|cc}
                $(N,s')$ & 800 point sample & 1000 point sample\\
                \toprule
                $(10,40)$ & \raisebox{-.5\height}{\includegraphics[width=7cm]{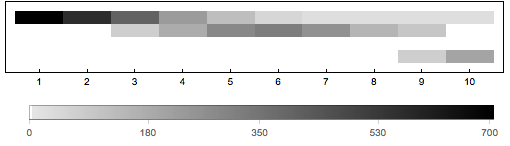}} & \raisebox{-.5\height}{\includegraphics[width=7cm]{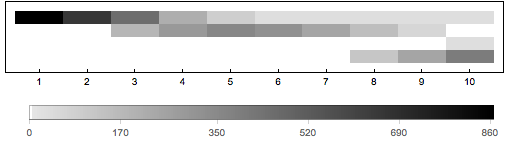}}\\[1.5cm]
                $(20,20)$ & \raisebox{-.5\height}{\includegraphics[width=7cm]{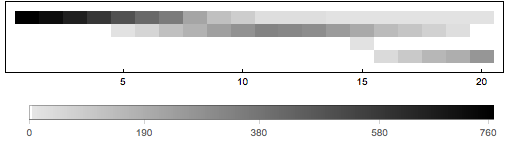}} & \raisebox{-.5\height}{\includegraphics[width=7cm]{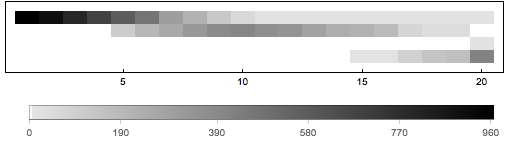}}\\[1.5cm]
                $(40,10)$ & \raisebox{-.5\height}{\includegraphics[width=7cm]{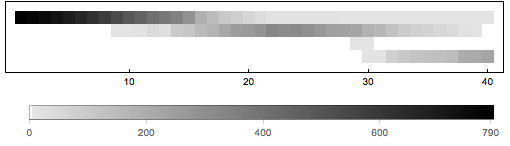}} & \raisebox{-.5\height}{\includegraphics[width=7cm]{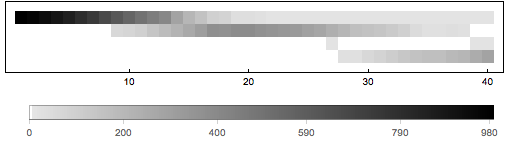}}\\[1.5cm]
                $(80,5)$ & \raisebox{-.5\height}{\includegraphics[width=7cm]{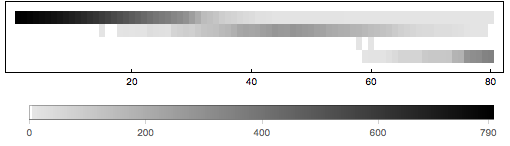}} & \raisebox{-.5\height}{\includegraphics[width=7cm]{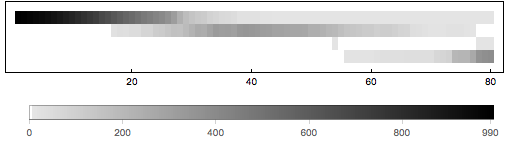}}
            \end{tabular}
            \caption{Here $N$ is the \emph{number of steps} and $s'=s\cdot10^4$ is the (scaled) step size. The top horizontal bar in each barcode represents $H_0$, and the bottom represents $H_3$.}\label{tb:cow-barcodes}
        \end{table}


    \section{Conclusions} 
    \label{sec:conclusions}

        \begin{quotation}
            \raggedleft
            \emph{``Life and death are one thread, the same line viewed from different sides.''}\\
            -- Lao Tzu, \emph{Tao Te Ching}
        \end{quotation}

        In this paper we have developed the tools of (extended) persistent homology using the category-theoretic language of \cite{Bubenik:dn}, formulated the relevant stability theorems, and presented various worked examples (focusing largely on the Klein bottle immersed\footnote{
            If we had had more time, it would have been an interesting project to study other immersions of the Klein bottle in $\rr^3$, such as the figure-8 immersion, for example.
        } in $\rr^3$) and commentary.
        There are some potentially interesting questions to consider in regards to the material that we have covered, though the relevance, the importance, and even the validity of them is, to the author, uncertain.

        \begin{enumerate}
            \item \emph{Distribution of point sampling:} We mentioned in \cref{sec:persistent_homology_of_a_cow} that, if we were handed a uniformly sampled point cloud from a space, then we probably wouldn't gain much from asking for another set of uniformly sampled data.
                It seems like there should be some way of analysing how the distribution of the sampled data  affects the resulting barcode: if we know that our points are distributed unevenly, say polynomially or exponentially more points are found around certain areas, can we predict how this will change the barcode?
                Obviously, for persistent homology, if we are dealing with some smooth manifold then it is the critical points that are of interest.
                So if we are told that our points were sampled from near these critical points, then how can we use this fact to improve our method of computing persistent homology?
                Depending on the specific embedding or immersion of the manifold, a naive calculation of the Vietoris-Rips complex might provide woefully inaccurate estimations as to the global homology.
            \item \emph{Symmetry of extended persistent homology:} The barcodes in \cref{fg:klein-bifiltration-barcode} have an interesting symmetry to them.
                Part of this is easily understandable: with this `convex' immersion of the Klein bottle and any height function, the $0$th homology will always be born at $-M$ and die at $M+1$, and the $2$nd homology will always be born at $M$ and die at $3M+1$.
                But how can we formalise the symmetry in the $1$st homology?
                If we compute the extended persistence homology of the torus embedded `vertically' in $\rr^3$ (i.e. with the hole perpendicular to the height function) then the intervals of the two $1$st homology classes have a rotational symmetry around $M+\frac12$, which is different to the mirror symmetry found in \cref{fg:klein-bifiltration-barcode}.
                Is this possibly to do with the orientability of the surface?
            \item \emph{Software for extended persistent homology:} As listed in \cite{Otter:2015wi}, there are many libraries for computing persistent homology, all with their various strengths and weaknesses.
                None of them, however, seem to be able to compute extended persistent homology.
                Would it be feasible to extend them to be able to do so, or would it be easier to write new dedicated code to do this task?
                It also seems possible that, with the powerful tools provided by Mathematica, we could compute the (extended) persistent homology of 3D models: we can define a height function $f$; calculate the `height slices' $F(a)$; discretise (triangulate) the resulting space; convert into simplicial data; and then calculate the homology computationally.
            \item \emph{Probabilistic persistent homology:} Extended persistent homology ensures that all classes die in finite time.
                This means that, at any given time, the persistence of a class is simply `how long it has left to live'.
                If we are handed some $k$-th homology class at time $a$ then we can calculate a probability that it is alive at time $a+p$ by simply looking at the proportion of $k$-th homology classes that persist until $a+p$.
                This naive idea becomes more interesting when we turn it around: given a slice of the barcode, say between times $t$ and $t'$, can we use the probabilities of a class persisting for time $p$ calculated from that slice to predict what might happen elsewhere in the barcode?
                That is, if we see many $k$-th homology classes being born and dying then it seems possible that this is a trait of the filtration -- it is `$k$-th homology noisy' -- and so we can take into account the fact that these classes have a high probability of dying in a short time when trying to reconstruct the rest of the barcode from finite data.
                Alternatively, is there anything to be gained from adopting an entirely stochastic viewpoint of persistent homology: drawing influence from the theory of statistical lifetime models?
                We don't allow ourselves access to all of the information about $X$, or its barcode.
                Instead, we can randomly sample various homology classes and look at how long they persists.
                In doing so, we can build up a statistical idea of the persistent homology, and try to fit it to some probabilistic model.
        \end{enumerate}



    \clearpage

    \begin{appendices}
        \crefalias{section}{appsec}
            \section{Cellular homology and Poincaré duality} 
    \label{sec:appendix_cells}
    
        \subsection{CW complexes and cellular homology} 
        \label{sub:cw_complexes_and_cellular_homology}
    
            \emph{Most of our definitions in this section come from \cite{hatcher2002algebraic}.}

            Here we briefly summarise the ideas behind \emph{CW complexes} and \emph{cellular homology}, and also agree on various notational quirks.
            The reader with a working knowledge of these topics can skip this section, referring back to it only when confused by potentially unfamiliar notation.

            \begin{definition}
                For $m\in\nn$ let $N_m=\{n\in\nn \mid n\leqslant m\}=\{0,1,\ldots,m\}$ and define $N_\infty=\nn$.
                Call a set of the form $N_k$ (for $k\in\nn\cup\{\infty\}$) a \emph{natural interval}\footnote{
                    This is not standard terminology, but we introduce it here to simplify certain statements throughout this paper.
                }.
            \end{definition}

            \begin{definition}[CW complex]
                Let $N_k$ be a natural interval, and let $\{e_\alpha^n\}_{\alpha\in \alpha_n}$ be a non-empty set of $n$-cells (copies of open $n$-discs $D^n$) for each $n\in N_k$.
                We build a topological space $X$, called a \emph{cell} (or \emph{CW}\footnote{
                    \emph{C} stands for `closure-finite' and \emph{W} stands for `weak topology'.
                }) \emph{complex}, by the following inductive procedure:
                \begin{enumerate}[(i)]
                    \item define $X^0=\{e_\alpha^0\}$;
                    \item define the \emph{$n$-skeleton $X^n$} by attaching\footnote{
                        i.e. take the quotient by an equivalence relation: $X^n=\big(X^{n-1}\sqcup_\alpha \overline{D_\alpha^n}\big)/\{x\sim\varphi_\alpha^n(x)\}$; the attaching map tells us how the boundary of the closed $n$-disc gets mapped into $X^{n-1}$.
                    } each $e_\alpha^n$ to $X^{n-1}$ via a map
                        \begin{equation*}
                            \varphi_\alpha^n\colon S^{n-1}\to X^{n-1};
                        \end{equation*}
                    \item If $k\in\nn$ then we set $X=X^k$; if $k=\infty$ then we set $X=\bigcup_n X^n$ and endow $X$ with the \emph{weak topology}: a set $U\subset X$ is open if and only if $U\cap X^n$ is open for all $n\in\nn$.
                \end{enumerate}
                Define the \emph{CW-dimension} of $X$ as \mbox{$\dim_\text{CW}(X)=k$}.
                We write $\sigma^n(X)$ to mean the underlying set structure of the $n$-skeleton:
                \begin{equation*}
                    \sigma^n(X) =
                    \begin{cases}
                        \bigcup_{m=0}^n\{e_\alpha^m\}_{\alpha\in\alpha_m} & n\leqslant\dim_\text{CW}(X);\\[.2cm]
                        \sigma^{\dim_\text{CW}(X)}(X) & n>\dim_\text{CW}(X),
                    \end{cases}
                \end{equation*}
                and define $\sigma(X)=\sigma^\infty(X)$.
            \end{definition}

            Given a polygon with side identifications we can realise it as a CW complex as follows.
            Say\footnote{
                This construction can be done in a more general case -- see \cite[Cell~Complexes,~Chapter~0,~p.~5]{hatcher2002algebraic} -- but all the examples that we encounter will be of this form.
            } the polygon has $2g$ edges, with sides identified in pairs, i.e. the boundary word is some permutation\footnote{
                e.g. if $2g=6$ then permissible boundary words include $aabbcc$ and $aba^{-1}cbc$.
            } of $a_1^{\pm1}a_1^{\pm1}\ldots a_g^{\pm1}a_g^{\pm1}$.
            Constructing the resulting surface is equivalent to attaching $g$ $1$-cells to a $0$-cell with the constant attaching map, which gives a wedge sum of $g$ copies of $S^1$.
            We label the $i$th copy of $S^1$ with $a_i$.
            Then we attach a $2$-cell to the wedge sum $\vee_{i=1}^g S^1$ along the boundary word\footnote{
                Split $S^1$ into $2g$ regions: $R_j=\left\{(\cos\theta,\sin\theta) \mid \theta\in\left[\frac{j\pi}{g},\frac{(j+1)\pi}{g}\right)\right\}$ for $j=0,\ldots,2g-1$.
                Say that the first letter of the boundary word is $a_i^{\vare_1}$, where $\vare_1=\pm1$, and pick some orientation for the copies of $S^1$ in the wedge sum $X^1=\vee_{m=1}^g S^1$.
                Then we define $\varphi^2\colon S^1\to X^1$ by mapping, from endpoint to endpoint, $R_0$ onto the copy of $S^1$ labelled with $a_i$, reversing the orientation if $\vare_1=-1$.
            }.

            \bigskip
            This procedure is hopefully made clear in the following example.

            \begin{figure}[t]
                \centering
                \includegraphics[width=.6\textwidth]{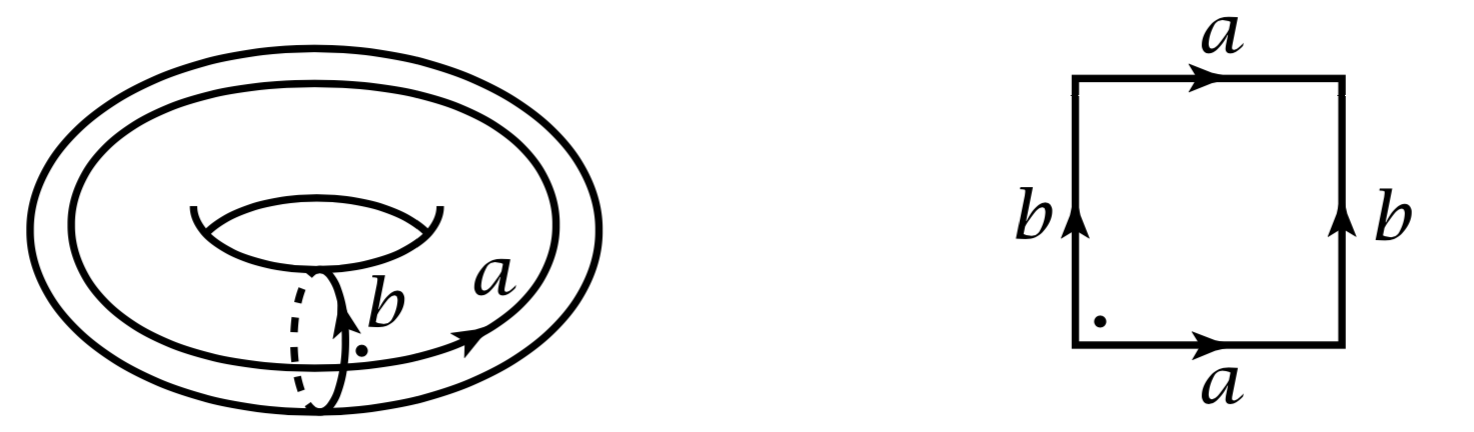}
                \caption{The torus as a polygon with side identifications \cite[p.~5,~\S0]{hatcher2002algebraic}.}\label{fg:hatcher-torus}
            \end{figure}

\begin{figure}[ht]
                \centering
                \frame{\includegraphics[width=.9\textwidth]{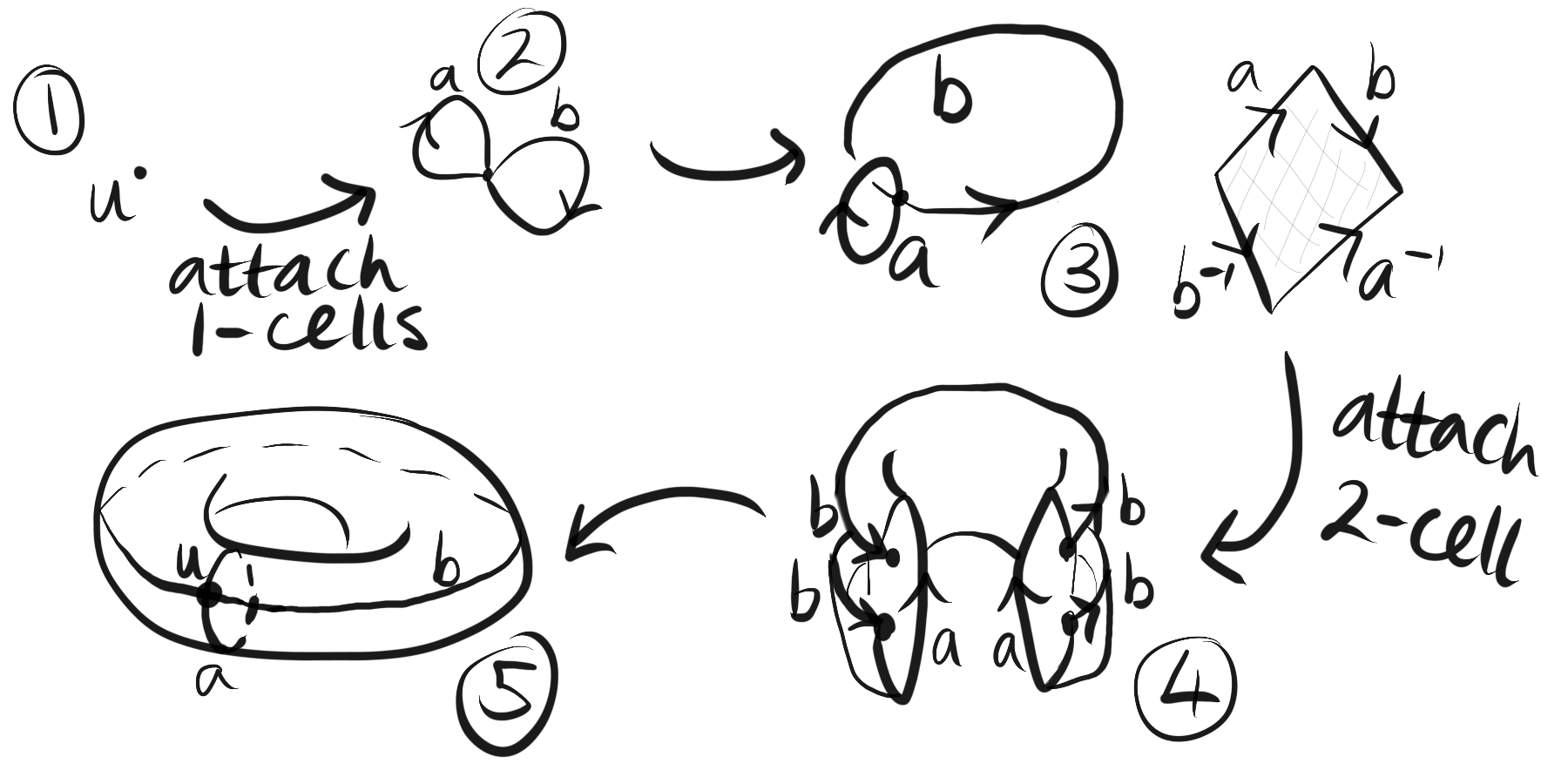}}
                \caption{Building the torus as a $CW$ complex, $n$-cell by $n$-cell.}\label{fg:tim-torus}
            \end{figure}

            \begin{example}[CW-complex structure of the torus]
                The torus $T$ can be defined as the surface resulting from the square with boundary word $aba^{-1}b^{-1}$ after side identification (see \cref{fg:hatcher-torus}).
                This induces a CW-complex structure on $T$, where
                \begin{itemize}
                    \item $\dim_\text{CW}(T) = 2$;
                    \item $\sigma(T) = \{e^0_0\}\cup\{e^1_0,e^1_1\}\cup\{e^2_0\}$;
                    \item $\varphi_\alpha^1\colon S^0=\{-1,1\}\to T^0=\{e_0^0\}$ is the constant map;
                    \item $\varphi_0^2\colon S^1\to T^1=\vee_{i=1}^2 S^1$ is the map $aba^{-1}b^{-1}$.
                \end{itemize}
                The map $\varphi_0^2$ is realised as follows (see \cref{fg:tim-torus}).
                We label one of the circles in the wedge sum $a$ and the other one $b$, label the point where they are identified $v$, and choose some orientation: say, clockwise.
                Next we cut $S^1$ at a point, resulting in a line, and identify (`glue') one end to $v$.
                Then we glue the first quarter of this line to the circle $a$ in a clockwise manner.
                We glue the next quarter to the circle $b$, also clockwise.
                The third quarter gets glued to the circle $a$ again, but this time anticlockwise.
                Finally, the last quarter gets glued to circle $b$, but again anticlockwise.
            \end{example}

            To define cellular homology of a CW complex, we first need to define the \emph{associated cellular chain complex}.

            \begin{definition}[Cellular chain complex]
                Let $X$ be a CW complex.
                Define the \emph{cellular chain complex $C_\bullet^{\text{CW}}(X)$ of $X$} as
                \begin{equation*}
                    C_\bullet^{\text{CW}}(X) = \ldots \xrightarrow{d_{n+1}} C_n^{\text{CW}}(X) \xrightarrow{d_n} C_{n-1}^{\text{CW}}(X) \xrightarrow{d_{n-1}} \ldots
                \end{equation*}
                where
                \begin{equation*}
                    C_n^\text{CW}(X) = H_n(X^n,X^{n-1})
                \end{equation*}
                and the $d_n$ are compositions coming from the long exact sequences for the pairs $(X^k,X^{k-1})$ and using \cite[Lemma~2.34,~\S2.2]{hatcher2002algebraic}: see \cite[\S2.2,~p.~139]{hatcher2002algebraic}.
            \end{definition}

            By \cite[Lemma~2.34,~\S2.2]{hatcher2002algebraic} we know that $C_n^\text{CW}(X)$ is free abelian and can be thought of as being generated by the $n$-cells of $X$.

            \begin{definition}[Cellular homology]
                Let $X$ be a CW complex.
                Define \emph{$n$-th cellular homology group $H_n^\text{CW}(X)$} as the $n$-th homology group of the associated cellular chain complex $C_\bullet^{\text{CW}}(X)$, i.e.
                \begin{equation*}
                    H_n^\text{CW}(X) = \frac{\ker d_n}{\im d_{n+1}}.\qedhere
                \end{equation*}
            \end{definition}

            \begin{theorem}\label{th:cellular-iso-sing-or-simp}
                Let $X$ be a CW complex.
                Then the $n$-th cellular homology group is isomorphic to the $n$-th homology\footnote{
                    Recall \cref{th:simplicial-iso-singular} -- we can use either simplicial or singular homology.
                } group, i.e.
                \begin{equation*}
                    H_n^\text{CW}(X)\cong H_n(X).\qedhere
                \end{equation*}
            \end{theorem}

            The above theorem means that, given some CW complex $X$, we can write $H_n(X)$ to mean \emph{the $n$-th homology group of $X$} without specifying whether we mean singular, simplicial, or cellular\footnote{
                Though sometimes, to clarify what we mean, we do use the notation $H_n^\text{CW}$.
            } -- they are all the same, up to isomorphism.
            Similarly, we can use any of these three types of homology to perform explicit calculations.

            \begin{definition}[Attach-and-collapse map]
                Let $X$ be a CW complex, $e_\alpha^n$ an $n$-cell, and $e_\beta^{n-1}$ an $n-1$-cell.
                Define the \emph{attach-and-collapse map}\footnote{
                    Again, this is not standard terminology.
                } \emph{$\chi^n_{\alpha\beta}\colon S_\alpha^{n-1}\to S_\beta^{n-1}$} as the composition $\chi^n_{\alpha\beta}=q_\beta^{n-1} q^{n-1}\varphi_\alpha^n$, where
                \begin{itemize}
                    \item $\varphi_\alpha^n$ is the attaching map;
                    \item $q^{n-1}\colon X^{n-1}\to X^{n-1}/X^{n-2}$ is the quotient map;
                    \item $q_\beta^{n-1}\colon X^{n-1}/X^{n-2}\to S_\beta^{n-1}$ is the map that collapses the complement of $e_\beta^{n-1}$ to a point\footnote{
                        For precise details, see \cite[141]{hatcher2002algebraic}.
                    }.\qedhere
                \end{itemize}
            \end{definition}

            \begin{theorem}[Cellular boundary formula]\label{th:cellular-boundary-formula}
                Let $X$ be a CW complex and
                \begin{equation*}
                    C_\bullet^{\text{CW}}(X) = \ldots \xrightarrow{d_{n+1}} C_n^{\text{CW}}(X) \xrightarrow{d_n} C_{n-1}^{\text{CW}}(X) \xrightarrow{d_{n-1}} \ldots
                \end{equation*}
                its associated cellular chain complex.
                Then, for $n>1$, the boundary maps are given by
                \begin{equation*}
                    d_n\colon e^n_\alpha \mapsto \sum_{\beta\in \alpha_{n-1}} \deg\left(\chi^n_{\alpha\beta}\right) e^{n-1}_\beta
                \end{equation*}
                and $d_1$ is the same as the simplicial boundary map $d_1^\Delta\colon\Delta_1(X)\to\Delta_0(X)$.
            \end{theorem}

            \begin{proof}
                See \cite[pp.~140,~141,~\S2.2]{hatcher2002algebraic} for the case where the homology coefficients are in $\zz$, and see \cite[Lemma~2.49,~\S2.2]{hatcher2002algebraic} for how to apply this with general group coefficients.
            \end{proof}




        \subsection{Orientability and Poincaré duality} 
        \label{sub:orientability_and_poincare_duality}

            Given some $n$-manifold we can generalise the idea of \emph{orientability} to \emph{$R$-orientability}, where $R$ is some commutative ring with identity, and in such a way that $\zz$-orientability recovers the original notion exactly.
            In particular, it can be shown that every manifold is $\zz_2$-orientable.
            See \cite[p.~235,~\S3.3]{hatcher2002algebraic} for explicit details.

            \begin{theorem}[Poincaré Duality]\label{th:poincare-duality}
                Let $M$ be an $R$-orientable closed $n$-manifold.
                Then\footnote{
                    The theorem actually hands us an explicit isomorphism between these homology groups: see \cite[Theorem~3.30,~\S3.3]{hatcher2002algebraic}.
                    We don't mention it here since we have no need for it.
                }
                \begin{equation*}
                    H^k(M;R)\cong H_{n-k}(M;R).\qedhere
                \end{equation*}
            \end{theorem}

            \begin{proof}
                See \cite[~pp.~247,~248,~\S3.3]{hatcher2002algebraic}.
            \end{proof}



    \end{appendices}


\addcontentsline{toc}{section}{References}

    \begin{thebibliography}{99}

        \bibitem{Bubenik:dn}
            Peter Bubenik and Jonathan A Scott.
            ``Categorification of persistent homology''.
            In: \emph{arXiv.org} (Jan. 2014).
            arXiv: \texttt{1205.3669v3 [math.AT]}.

        \bibitem{CohenSteiner:2007is}
            David Cohen-Steiner, Herbert Edelsbrunner, and John Harer.
            ``Stability of Persistence Diagrams''.
            In: \emph{Discrete \& Computational Geometry} 37.1 (2007), pp.~103-120.

        \bibitem{CohenSteiner:2009ho}
            David Cohen-Steiner, Herbert Edelsbrunner, and John Harer.
            ``Extending Persistence Using Poincaré and Lefschetz Duality''.
            In: \emph{Foundations of Computational Mathematics} 9.1 (2009), pp.~79-103.

        \bibitem{Edelsbrunner:2008gf}
            Herbert Edelsbrunner and John Harer.
            ``Persistent homology – a survey''.
            In: \emph{Surveys on Discrete and Computational Geometry}.
            Ed. by János Pach, Jacob E Goodman, and Richard Pollack.
            American Mathematical Society, 2008, pp.~257-282.
            {\sc ISBN:} 978-0-8218-4239-3.

        \bibitem{ghrist2014elementary}
            Robert Ghrist.
            \emph{Elementary Applied Topology}.
            Createspace, 2014.
            {\sc ISBN:} 978-1-5028-8085-7.

        \bibitem{hatcher2002algebraic}
            Allen Hatcher.
            \emph{Algebraic Topology}.
            Cambridge University Press, 2002.
            {\sc ISBN:} 978-0-521-79540-1.

        \bibitem{Massey:1967we}
            William S Massey.
            \emph{Algebraic topology: an introduction}.
            Harcourt, Brace \& World, 1967.

        \bibitem{Perseus:70k_p2_v}
            Vidit Nanda.
            \emph{Perseus, the Persistent Homology Software}.
            {\sc URL:} \url{http://www.sas.upenn.edu/~vnanda/perseus}.

        \bibitem{Neeb:2004te}
            Karl-Hermann Neeb.
            ``Current groups for non-compact manifolds and their central extensions''.
            In: \emph{Infinite Dimensional Groups and Manifolds}.
            Ed. by Tilman Wurzbacher.
            2004, pp.~109-184.
            {\sc ISBN:} 978-3-11-018186-9.

        \bibitem{Otter:2015wi}
            Nina Otter, Mason A Porter, Ulrike Tillmann, Peter Grindrod, and Heather A Harrington.
            ``A roadmap for the computation of persistent homology''.
            In: \emph{arXiv.org} (June 2015).
            arXiv: \texttt{1506.08903v3 [math.AT]}.

        \bibitem{Mathematica:tp}
            Wolfram Research, Inc.
            \emph{Mathematica}.
            Version 10.4.

    \end{thebibliography}
\end{document}